\documentclass[final]{amsart}

\usepackage{booktabs}
\usepackage{cite}

\usepackage{lipsum}
\usepackage{amsfonts}
\usepackage{graphicx}
\usepackage{epstopdf}
\usepackage{algorithmic}
\usepackage{algorithm}
\usepackage{amsopn}
\usepackage{amssymb}
\usepackage{amsmath}
\usepackage{amsthm}
\usepackage[colorlinks=true]{hyperref}
\usepackage{enumerate}
\usepackage{todonotes}
\usepackage{tikz-cd}
\usepackage{diagbox}
\usepackage{multirow}
\ifpdf
  \DeclareGraphicsExtensions{.eps,.pdf,.png,.jpg}
\else
  \DeclareGraphicsExtensions{.eps}
\fi

\newtheorem{theorem}{Theorem}[section]

\newtheorem{corollary}[theorem]{Corollary}
\theoremstyle{definition}
\newtheorem{definition}[theorem]{Definition}

\theoremstyle{remark}
\newtheorem{remark}[theorem]{Remark}

\def\th{\theta}

\def\la{\lambda} 
\def\rh{\rho}

\def\ph{\varphi}

\def\SE{S\hspace{-0.08em}E}

\def\x{\times}

\def\L{\mathcal{L}}

\let\on=\operatorname

\def\R{\mathbb{R}}

\let\on=\operatorname

\def\Imm{\on{Imm}}
\newcommand{\ud}{\,\mathrm{d}}
\newcommand{\ip}[2]{\langle #1,#2 \rangle}

\begin{document}

\title{A relaxed approach for curve matching with elastic metrics}

\author{Martin Bauer}
\address{Florida State University}
\email{bauer@math.fsu.edu}

\author{Martins Bruveris}
\address{Brunel University London}
\email{martins.bruveris@brunel.ac.uk}

\author{Nicolas Charon}
\address{Johns Hopkins University}
\email{charon@cis.jhu.edu}

\author{Jakob M\O ller-Andersen}
\address{Florida State University}
\email{jmoeller@math.fsu.edu}

\renewcommand{\shortauthors}{M. Bauer, M. Bruveris, N. Charon and J. M\O ller-Andersen}



\maketitle

\begin{abstract}
In this paper we study a class of Riemannian metrics on the space of unparametrized curves and develop a method to compute geodesics with given boundary conditions. It extends previous works on this topic in several important ways. The model and resulting matching algorithm integrate within one common setting both the family of $H^2$-metrics with constant coefficients and scale-invariant $H^2$-metrics on both open and closed immersed curves. These families include as particular cases the class of first-order elastic metrics. An essential difference with prior approaches is the way that boundary constraints are dealt with. By leveraging varifold-based similarity metrics we propose a relaxed variational formulation for the matching problem that avoids the necessity of optimizing over the reparametrization group. Furthermore, we show that we can also quotient out finite-dimensional similarity groups such as translation, rotation and scaling groups. The different properties and advantages are illustrated through numerical examples in which we also provide a comparison with related diffeomorphic methods used in shape registration.
\end{abstract}

\section{Introduction}
In this article we study Riemannian metrics on the  space of unparametrized, $\R^d$-valued curves. The interest in this topic is fueled by applications in medical imaging, computer animation, geometric morphometry and other fields~\cite{Srivastava2016}. The space of closed curves is important in shape analysis where it is used to study objects that can be represented by the shape of their boundary \cite{Younes1998,Klassen2004,Laga2014,BBHM2015a}. At the same time open curves are relevant in applications such as the analysis of hurricane paths, bird migration patterns \cite{su2014,SBK2017}, and human character motions \cite{Esl2014,Esl2014b} or in character and speech recognition \cite{su2014b}.

The analysis of shapes and their differences relies on the notion of a distance. To define such a distance, we will start with a Riemannian metric on the space of curves and use the induced geodesic distance to quantify differences between curves.
Mathematically we model curves as smooth mappings from a parameter space to $\R^d$. The parameter space is $S^1$ for closed curves and $[0,2\pi]$ for open curves. On the space of curves we can consider the action of the reparametrization group and we consider two curves equivalent if they only differ by a reparametrization, i.e., $c_1\sim c_2$ if $c_1=c_2\circ\varphi$ for some reparametrization $\varphi$. To define a Riemannian metric on the space of unparametrized curves we will start with a reparametrization invariant metric on the space of parametrized curves and consider the induced metric on the quotient space of unparametrized curves\footnote{Note, that the invariance of the metric is only a necessary condition for a metric on the space of parametrized curves to induce a metric on unparametrized curves. However, all the metrics considered in this article do induce Riemannian metrics on the quotient space. For details, see \cite{Michor2007}.}. 

The simplest invariant metric on the space of parametrized curves is the $L^2$-metric
\begin{equation*}
G_c(h,k)=\int_{M^1} \langle h,k \rangle \ud s \;;
\end{equation*}
where $c$ is a curve, $h,k$ are tangent vectors to $c$, $\langle\cdot,\cdot\rangle$ is the Euclidean inner product, $\ud s = |c'| \ud \theta$ is integration with respect to arc length and $M^1$ is either $S^1$ or $[0,2\pi]$.
It came as a big surprise when Michor and Mumford found in \cite{Michor2006c} that this metric induces vanishing geodesic distance on the space of unparametrized curves\footnote{This result has been later extended to the space of parametrized curves in \cite{Bauer2012c}.}. Here vanishing geodesic distance means that given any two curves there exist paths of arbitrary short length connecting them and consequently the geodesic distance, which is defined as the infimum over all path lengths, is identically zero. This result renders the $L^2$-metric impractical for applications in shape analysis and thus started the quest for stronger and more meaningful metrics:
Michor and Mumford proposed curvature weighted versions of the $L^2$-metric \cite{Michor2006c}  and Shah \cite{Shah2008} studied length weighted versions to successfully overcome the degeneracy of vanishing geodesic distance. A more promising approach, first investigated out of purely theoretical interest, is to include derivatives of the tangent vector in the Riemannian metric, yielding the class of Sobolev metrics \cite{Michor2007,Mennucci2008}
\begin{equation*}
G_c(h,k)=\int_{M^1} \langle h,k \rangle \ud s+\langle D_s^n h,D_s^n k \rangle \ud s\;,\qquad n\geq 1\;,
\end{equation*}
where $D_s=\frac{1}{|c'|}\partial_{\theta}$ denotes derivative with respect to arc length. 
Closely related to first order Sobolev metrics is the family of elastic $G^{a,b}$-metrics, proposed as a model for shape analysis in \cite{Mio2004,Mio2007}. These are metrics of the form
\begin{equation*}
G_c(h,k) = 
\int_{M^1}
a \langle D_s h^\top,D_s k^\top \rangle 
+
b \langle D_s h^{\bot},D_s k^{\bot} \rangle  \ud s \,.
\end{equation*}
Here $a, b > 0$ are constants and $D_s h^\top = \langle D_s h, D_s c\rangle D_s c$ and $D_s h^{\bot} = D_s h - D_s h^\top$ denote the decomposition of $D_s h$ into vectors tangent and orthogonal to the curve $c$. The first term can be interpreted as penalizing stretching of the curve, whereas the second term measures bending of the curve. Thus one expects to be able to model a variety of behaviors by varying the constants $a$ and $b$. However, so far only two special cases have been implemented numerically\footnote{The article \cite{Mio2009} studied the family of elastic metrics with different choices of parameters on the space of \emph{parametrized} curves. This analysis was, however, not extended to the space of \emph{unparametrized} curves, which is the more relevant object for applications.}: Younes et. al.~\cite{Michor2008a} studied the metric for planar curves when $a=b=1$ and Srivastava et al.~\cite{Jermyn2011,Klassen2004} studied it for $\R^d$-valued curves when $a=1$ and $b=1/4$. In both cases there exist transformations that allow one to compute explicit formulas for geodesics and geodesic distance. These transformations allowed for the development of efficient and fast numerical algorithms both on parametrized and unparametrized curves and have been proven successful for applications in shape analysis. 
In \cite{Bauer2014b} these transformations have been extended to a wider class of parameters and metrics. In an upcoming preprint Kurtek and Needham \cite{KuNe2018} propose a different numerical framework for the class of $G^{a,b}$-metrics based on a generalization of the transformation of Younes et. al. \cite{Michor2008a}.

It has been recently shown in \cite{Bruveris2014b_preprint,Bruveris2014,Vialard2014_preprint}, that adding second derivatives to the metric allows one to obtain completeness results for the Riemannian manifolds in question: the geodesic equation is globally well-posed, the metric completion consists of all $H^2$-immersions, and the metric extends to a strong Riemannian metric on the space of $H^2$-immersions. Furthermore, any two curves in the same connected component can be joined by a minimizing geodesic. For a more detailed overview of various Riemannian metrics on the space of curves we refer to the overview articles \cite{Bauer2014,BBM2016}. 

\subsection*{Contributions of the article}
In this article we present the first numerical implementation of the geodesic initial and boundary value problem for a family of first and second order metrics on the space of open and closed unparametrized curves. Our code is available under an open source license\footnote{\url{https://www.github.com/h2metrics/h2metrics}}. This family includes in particular the elastic $G^{a,b}$-metrics with arbitrary parameters $a$ and $b$ as well as scale-invariant Sobolev metrics. To be precise, we study metrics of the form
\begin{multline*}
G_c(h,k) = 
\int_{M^1} a_0(\ell_c) \langle h,k \rangle+
a_1(\ell_c)\langle D_s h^\top,D_s k^\top \rangle 
+
b_1(\ell_c)\langle D_s h^{\bot},D_s k^{\bot} \rangle 
\\+a_2(\ell_c)\langle D_s^2 h,D_s^2 h \rangle \ud s \,.
\end{multline*}
where $a_0, a_1, b_1, a_2 \in C^{\infty}(\R_{>0},\R_{\geq 0})$ are smooth positive functions of the curve length $\ell_c$. All metrics in this class are invariant with respect to Euclidean motions thus they induce Riemannian metrics on the shape space of unparametrized curves modulo Euclidean motions. If the coefficient functions $a_j$ and $b_1$ are chosen appropriately one obtains scale-invariant metrics, which then induce Riemannian metrics on the space of  unparametrized curves modulo similarity transformations.

In future applications this numerical framework will allow us to choose the constants of the metric in a data-driven way as opposed to the ad hoc methods employed currently. As a first step towards this goal we will show in selected examples how the choice of constants influences minimal geodesics between two shapes and how it affects the point-to-point registration. In the experiments section we also compare this intrinsic approach to curve matching with the LDDMM framework, where differences between curves are measured extrinsically using a Riemannian metric on the diffeomorphism group of the ambient space.  

On the space of closed curves we additionally extend the completeness results, that were obtained first in \cite{Bruveris2014, Bruveris2014b_preprint} for Sobolev metrics with constant coefficients and then in \cite{Bruveris2017_preprint} for length-weighted Sobolev metrics to this class of elastic metrics. For open curves we find a counter example showing that second order metrics with constant coefficients are not metrically complete.

From a numerical point of view, we introduce a different method for handling the boundary conditions when solving the geodesic boundary value problem on the space of unparametrized curves. Mathematically, an unparametrized curve corresponds to the orbit $c \circ \on{Diff}(M^1)$ but there is no numerically convenient way to discretize this group action. This problem has been approached in various ways, e.g. in previous work \cite{BBHM2017} the authors optimized simultaneously over both the geodesic path and the reparametrization of the target curve. In this paper we continue to develop the idea inspired from previous works on curve and surface registration based on diffeomorphic models like \cite{Glaunes2006,Glaunes2008,Charon2013}: to enforce the constraint that the endpoint of the geodesic path $c(1)$ and the given target curve $c_1$ belong to the same equivalence class of unparametrized curves, we use an auxiliary reparametrization-invariant distance function. The construction of those distance functions, unrelated to Sobolev metrics, follow the principles of geometric measure theory, which have so far been used as fidelity terms in combination with models like the Large Deformation Diffeomorphic Metric Mapping (LDDMM) framework. In the present article, building on our previous conference publication \cite{BauerMICCAI2017}, we examine a fairly general class of kernel metrics on immersed open and closed curves that are induced from the representation of curves as oriented varifolds. We examine the rigorous conditions to obtain distance functions on unparametrized immersed curves. Since this framework provides us with smooth proximity measures between unparametrized curves, it can be used to formulate the geodesic boundary value problem without having to explicitly estimate the reparametrizations. 


Our numerical implementation takes full advantage of the flexibility provided by the varifold-based distance. We formulate two algorithms: an exact matching algorithm using an augmented Lagrangian approach and an inexact matching algorithm based on optimizing a relaxed functional that incorporates the geodesic energy and the varifold-based distance. The latter algorithm is more flexible when the given data is noisy and similar approaches are used in most deformation-based shape analysis frameworks.

\section{Metrics on spaces of open and closed curves}

\subsection{Shape spaces of curves}\label{sec:shape_space}
Let $d \geq 2$ be the dimension of the ambient space and $M^1$ the parameter space for a curve. For open curves we have $M^1=I=[0,2\pi]$ and for closed curves $M^1 = S^1 = \R/2\pi\mathbb Z$. In both cases $M^1$ is a compact, one-dimensional manifold.
For a curve $c: M^1\to \mathbb R^d$ we write $c'(\theta)=\frac{d}{d\theta}c(\theta)$ to denote its derivative.
\begin{definition}
Let $M^1$ be $S^1$ or $I$. The space of smooth, regular curves with values in $\R^d$ is 
\begin{align*}
\Imm(M^1,\R^d)=\left\{c\in C^{\infty}(M^1,\R^d)\colon  c'(\th) \neq 0\quad \forall \th \in M^1 \right\}\,.
\end{align*}
\end{definition}
The notation $\Imm$ stands for \emph{immersions}. The space $\Imm(M^1,\R^d)$ is an open subset of the Fr\'echet space $C^\infty(M^1,\R^d)$ and therefore itself a Fr\'echet manifold. Its tangent space $T_c\Imm(M^1,\R^d)$ at any curve $c$ is the vector space $C^\infty(M^1,\R^d)$ itself.

We will call  curves in $\Imm(M^1,\R^d)$ parametrized curves, because as maps from the parameter space $M^1$ to $\mathbb R^d$ they carry with them a parametrization. We will later define the space of unparametrized curves in Definition~\ref{def:Bif_def}. 

Two curves that differ only by their parametrization represent the same geometric object. In the context of shape analysis it is therefore natural to consider them as equal, i.e., we identify the curves $c$ and $c \circ \ph$, where $\ph$ is a reparametrization. As the reparametrization group we use the group of smooth diffeomorphisms of $M^1$,
\begin{equation*}
\on{Diff}(M^1) = \left\{ \ph \in C^\infty(M^1,M^1) \,:\, \ph \text{ bij. and }
\ph^{-1}\in C^\infty(M^1,M^1)
\right\}\,,
\end{equation*}
which is an infinite-dimensional regular Fr\'echet Lie group \cite{Michor1997}.
For the two cases studied in this article these groups are
\begin{align*}
\on{Diff}(S^1) &= \left\{ \ph \in C^\infty(S^1,S^1) \,:\, \ph'(\th) \neq 0 \quad\forall \th \in S^1 \right\}\,,\\
\on{Diff}(I) &= \left\{ \ph \in C^\infty(I,I) \,:\, \ph'(\th) \neq 0,\; \ph(\{0,2\pi\})=\{0,2\pi\} \right\}\;.
\end{align*}
 
To define the quotient space of unparametrized curves we need to restrict ourselves to \emph{free immersions}, i.e., those upon which the diffeomorphism group acts freely:
\begin{definition}
Let $M^1$ be $S^1$ or $I$. The space of free immersions with values in $\mathbb R^d$ is 
\begin{equation*}
\Imm_f(M^1,\R^d)=\left\{c\in \Imm(M^1,\R^d):
\big(c \circ \ph = c \;\Rightarrow\; \ph = \on{Id}_{M^1}\big)
\right\}\,.    
\end{equation*}
\end{definition}

This restriction is only necessary for technical reasons to be able to define a manifold structure on the quotient space; in applications almost all curves are freely immersed, in particular the subset of free immersions is dense \cite{Michor1991}. 
\begin{definition}
\label{def:Bif_def}
The  space of unparametrized curves 
\begin{equation*}
 B_{i,f}(M^1,\R^d)=\Imm_f(M^1,\R^d) / \on{Diff}(M^1)\,,
\end{equation*}
is the orbit space of the group action of $\on{Diff}(M^1)$ restricted to all free immersions.
\end{definition}
This space is a Fr\'echet manifold although constructing charts is nontrivial in this case \cite{Michor1991}. The space  $B_{i,f}(M^1,\R^d)$ is sometimes referred to as the \emph{pre-shape space}, while the \emph{shape space} is obtained from the pre-shape space by additionally taking the quotient with respect to the group $S(d)$ of similarity transformations of $\R^d$ or one of its subgroups.
Here
\[
S(d) = \big(\R_{>0} \x SO(d)\big) \ltimes \R^d\,,
\]
where $\R_{>0}$ is the scaling group, $SO(d)$ is the rotation group and $\R^d$ is the translation group. The composition of two transformations is given by
\[
(r_1, A_1, w_1) \cdot (r_2, A_2, w_2) 
= (r_1r_2, A_1 A_2, r_2^{-1} A_2^{-1} w_1 + w_2)\,.
\]
The group $S(d)$ acts on curves from the left via
\begin{equation*}
\Imm(M^1,\R^d)\times S(d)\to \Imm(M^1,\R^d)\quad \left(c,(r,A,w)\right)\mapsto r A(c+w)\,.
\end{equation*}
Note that elements of $S(d)$ are orientation-preserving. We do not include the reflection $c \mapsto -c$ in this group. This can be done, but is not relevant for the applications considered below.

Let $H$ be a subgroup of $S(d)$.
Common choices for $H$ are the translation group $\R^d$, the group of Euclidean motions $SE(d) = SO(d) \ltimes \R^d$ and $S(d)$ itself. The  shape space of unparametrized curves modulo similarities of type $H$ is the quotient
\begin{equation*}
\mathcal S_H(M^1,\R^d)= B_{i,f}(M^1,\R^d)/H=\Imm_f(M^1,\R^d) / H \times \on{Diff}(M^1)\,.
\end{equation*}
We will write simply $\mathcal S(M^1,\R^d)$ instead of $\mathcal S_H(M^1,\R^d)$, when the meaning of $H$ is clear from the context.

\subsection{Notation}
We denote the Euclidean inner product on $\R^d$ by $\langle\cdot,\cdot\rangle$. For any fixed immersed curve $c$, we denote differentiation and integration with respect to arc length by $D_s=\frac 1{|c_\theta|}\partial_{\theta}$ and $\mathrm ds=|c_\theta|\ud \theta$ respectively. The length of the curve is $\ell_c=\int_{M^1}\ud s = \int_{M^1} |c'| \ud \th$. We will omit the subscript and write $\ell = \ell_c$ if the curve $c$ is clear from the context. The unit length tangent vector to $c$ is $v =D_s c = \frac{c'}{|c'|}$. We can decompose any vector field along the curve, $h \in T_c \on{Imm}(M^1,\R^d)$, into components tangential and normal to the curve,
$h = \langle h, v\rangle v + (h - \langle h, v \rangle v)$ and we denote them by $h^\top=\langle h,v \rangle v$ and $h^\bot = h-h^\top$. In fact we will apply this decomposition to the derivative $D_s h$ and we write $D_s h^\top = (D_s h)^\top$ and $D_s h^\bot = (D_s h)^\bot$. Note that $D_s$ does not commute with these projections and $D_s(h^\top) \neq (D_s h)^\top$ in general.

\subsection{Higher order elastic metrics}\label{sec:elastic_metrics_param}
Here we introduce the class of Riemannian metrics that will be used in the remainder of the article. We also show that these metrics possess certain desirable completeness properties and these will serve as theoretical justification for our numerical framework.

\begin{definition}
\label{def:sobolev_metric}
A \emph{second order elastic Sobolev metric with length-weighted coefficients} is a Riemannian metric on the space $\Imm(M^1,\R^d)$ of parametrized curves of the form
\begin{multline*}
G_c(h,k) = 
\int_{M^1} a_0(\ell) \langle h,k \rangle+
a_1(\ell)\langle D_s h^\top,D_s k^\top \rangle 
+
b_1(\ell)\langle D_s h^{\bot},D_s k^{\bot} \rangle 
\\+a_2(\ell)\langle D_s^2 h,D_s^2 h \rangle \ud s \,,
\end{multline*}
where $a_0, a_1, b_1, a_2 \in C^{\infty}(\R_{>0},\R_{\geq 0})$ are smooth positive functions of the curve length $\ell$ with $a_0(\ell) > 0$, $a_2(\ell) > 0$ and $h,k \in T_c\Imm(S^1,\R^d)$ are tangent vectors at $c$.
\end{definition}
In the remainder of this article we will restrict our attention to two special subfamilies. First, the family of \emph{elastic metrics with constant coefficients},
\begin{multline}
\label{eq:sobolev_metric:constant_coeff}
G^1_c(h,k) = 
\int_{M^1} a_0 \langle h,k \rangle+
a_1\langle D_s h^\top,D_s k^\top \rangle 
+
b_1\langle D_s h^{\bot},D_s k^{\bot} \rangle 
\\+a_2\langle D_s^2 h,D_s^2 h \rangle \ud s \;,
\end{multline}
where $a_0,a_1,b_1,a_2 \in \R_{\geq 0}$ and $a_0,a_2 > 0$;
second, the family of \emph{scale-invariant elastic metrics},
\begin{multline}\label{eq:sobolev_metric:scale_inv}
G^2_c(h,k) = 
\int_{M^1} \frac{a_0}{\ell^3} \langle h,k \rangle+
\frac{a_1}{\ell}\langle D_s h^\top,D_s k^\top \rangle 
+
\frac{b_1}{\ell}\langle D_s h^{\bot},D_s k^{\bot} \rangle 
\\+\ell a_2\langle D_s^2 h,D_s^2 h \rangle \ud s \;,
\end{multline}
where the coefficient functions are of the form $\ell \mapsto \la \ell^{\rh}$ and $a_0,a_1,b_1,a_2 \in \R_{\geq 0}$ with $a_0,a_2 > 0$. Note that the symbols $D_s$, $\ud s$, $\bot$ and $\top$ hide the nonlinear dependency of the metric on the base point $c$.
In  the following remark we discuss the invariance properties of the metrics 
\eqref{eq:sobolev_metric:constant_coeff} and \eqref{eq:sobolev_metric:scale_inv}.

\begin{remark}
 Because we use arc length operations in the definition of length-weighted elastic Sobolev metrics, Definition~\ref{def:sobolev_metric}, the resulting metrics are invariant under the action of the diffeomorphism group $\on{Diff}(M^1)$. They are also invariant under the Euclidean motion group $SE(d)$, but only the family $G^2$ of scale-invariant metrics is also invariant under the scaling group.
The invariance properties of these metrics will allow us later to define induced Riemannian metrics on the shape space of unparametrized curves. 
\end{remark}

First, however, we will study properties of length-weighted elastic metrics on the space of parametrized curves. Let $G$ be such a metric. The Riemannian length of a path $c\colon[0,1]\to\Imm(M^1,\R^d)$ is
\begin{align*}
L(c) = \int_0^1 \sqrt{G_{c(t)}(c_t(t),c_t(t))} \ud t\,,
\end{align*}
with $c_t=\partial_t c$ the time derivative of the path $c$. We denote by $\mathcal P$ the space of all smooth paths,
\[
\mathcal P = C^\infty([0,1], \on{Imm}(M^1,\R^d))\,.
\]
The geodesic distance induced by $G$ between two given curves $c_0$, $c_1$ is defined as the infimum of the lengths of all paths connecting these two curves, i.e., 
\[
\on{dist}^G(c_0, c_1) = \inf \left\{L(c)\,:\, c \in \mathcal P,\, c(0)=c_0,\, c(1)=c_1\right\}\,.
\]
It is a general result in Riemannian geometry that the squared geodesic distance is also the infimum over all paths of the Riemannian energy,
\begin{equation} \label{eq:EnergyFunctional}
\begin{aligned}
E(c) &= \int_0^1 G_{c(t)}(c_t(t), c_t(t)) \ud t \,.
\end{aligned}
\end{equation}
Geodesics are critical points of the energy functional and the first order condition for critical points, $DE(c) = 0$ is the geodesic equation.
For elastic metrics the geodesic equation is a partial differential equation for the function $c = c(t,\theta)$. 
Since we are working in infinite dimensions the existence of geodesics is a nontrivial question. For elastic Sobolev metrics we have the following existence results for geodesics, which are based on the results in \cite{Bruveris2014, Bruveris2014b_preprint, Michor2007}. They will serve as the theoretical foundation of the proposed numerical framework.
\begin{theorem}\label{thm:completeness}
Let $G$ be a second order elastic metric, either scale-invariant or with constant coefficients on the space of closed curves. Then
\begin{enumerate}
    \item The Riemannian
manifold $\left(\Imm(S^1,\R^d),G\right)$ is geodesically complete, i.e.,
given any initial conditions $(c_0,u_0)\in T \Imm(S^1,\R^d)$ the solution of the geodesic 
equation for the metric $G$ with initial values $(c_0,u_0)$
exists for all times.
   \item  The metric completion of the space $\Imm(S^1,\R^d)$ equipped with the geodesic distance 
   $\operatorname{dist}^{G}$ is the space $\mathcal I^2(S^1,\R^d)$ of immersions of Sobolev class $H^2$,
   \begin{equation*}
       \mathcal I^2(S^1,\R^d)=\left\{c\in H^{2}(S^1,\R^d)\,:\,  c'(\th) \neq 0\quad \forall \th \in S^1 \right\}\,.
   \end{equation*}
   Furthermore, any two curves in the same connected component of the space $\mathcal I^2(S^1,\R^d)$ can be joined by a minimizing geodesic.
\end{enumerate}
\end{theorem}

\begin{remark}[Incompleteness for open curves]
In the following we will present a counterexample for the above completeness result for the metric $G^1$ in the case of open curves. Therefore we consider the path 
$c(t,\theta)=((1-t)\theta,0)$ for $t\in[0,1]$. Note that this path leaves the space of immersions for $t=1$, since 
$c(1,\theta)=(0,0)$.  We have
\begin{equation*}
c_t = (-\theta,0),\quad  c_{\theta} = (1-t,0),\quad 
\langle D_s c_t,v\rangle = \frac{-1}{1-t},\quad  (D_s c_t)^{\bot} = (0,0)\quad D^2_s c_t =(0,0)\;.
\end{equation*}
Using this we calculate the $G^1$-length of $c$:
\begin{align*}
L(c) &= \int_0^1 \sqrt{G^1_{c(t)}(c_t(t),c_t(t))} \ud t\\&=
\int_0^1 \sqrt{\int_0^{2\pi} a_0\;\theta^2+\frac{a_1}{(1-t)^2}(1-t)\;\ud  \theta } \ud t\\
&=\int_0^1 \sqrt{ \;\frac{a_0(2\pi)^3(1-t)}{3}+\frac{a_12\pi}{(1-t)} } \ud t\\
&\leq C
\int_0^1 \sqrt{ 1+\frac{1}{1-t}} \ud t\\
&= C(\sqrt{2}+\operatorname{cosh}^{-1}(\sqrt{2}))
<\infty\,,
\end{align*}
where $C>0$ is a constant.
Thus we have found a path that leaves the space of immersions with finite $G^1$-length. This shows that the space
$\left(I^2([0,2\pi],\R^d), \operatorname{dist}^{G^1} \right)$ is metrically incomplete. 

We conjecture that for the scale invariant metric $G^2$ the completeness results would also hold on the space of open curves if  $a_0, a_1$ and $a_2$ are non-zero. Note that the above path $c$ has indeed infinite length with respect to the metric $G^2$.   
\end{remark}

\begin{proof}[Proof of Theorem~\ref{thm:completeness}]
Denote by $G^1$ an elastic metric with constant coefficients and by $G^2$ a scale-invariant elastic metric of the form~\eqref{eq:sobolev_metric:constant_coeff} and~\eqref{eq:sobolev_metric:scale_inv}. We first observe that both $G^1$ and $G^2$ extend to smooth Riemannian metrics on the Hilbert manifold $\mathcal I^2(S^1,\R^d)$. This follows directly from Sobolev embedding and multiplication theorems.

We will next show that $\mathcal I^2(S^1,\R^d)$ equipped with these metrics is metrically complete. For closed curves, metric completeness of second order Sobolev metrics on the space $\mathcal I^2(S^1,\R^d)$ has been shown for the metric
\begin{align*}
  \bar G^1_c(h,k) &= 
\int_{S^1}  \langle h,k \rangle+ \langle D_s h,D_s h \rangle + \langle D_s^2 h,D_s^2 h \rangle \ud s\,,
\end{align*}
in \cite{BH2015, Bruveris2014, Bruveris2014b_preprint} and for the metric
\begin{align*}
\bar G^2_c(h,k) &= 
\int_{S^1} \frac{1}{\ell^3} \langle h,k \rangle+\frac{1}{\ell}\langle D_s h,D_s h \rangle \ud s+\ell \langle D_s^2 h,D_s^2 h \rangle \ud s \,,
\end{align*}
in \cite{Bruveris2017_preprint}. Because we can find uniform constants $C_1$, $C_2$ such that
\begin{align*}
C_1 \bar G^1_c(h,h) \leq  G^1_c(h,h) \leq  C_2 \bar G^1_c(h,h)   \text{ and }C_1 \bar G^2_c(h,h) \leq  G^2_c(h,h) \leq  C_2\bar G^2_c(h,h)\;,
\end{align*}
holds for all $(c,h) \in T\mathcal I^2(S^1,\R^d)$, where $\bar G^1_c$, $\bar G^2_c$ are defined as above, it follows that the metrics $G^1$ and $G^2$ are also metrically complete on $\mathcal I^2(S^1,\R^d)$. 

Metric completeness of $G^1$ and $G^2$ on $\mathcal I^2(S^1,\R^d)$ implies geodesic completeness. To see that the space $\on{Imm}(S^1,\R^d)$ of smooth immersions is also geodesically complete, we use that the geodesic equation preserves smoothness of the initial conditions: if the initial curve and velocity field are $C^\infty$-smooth, then so is every curve along the geodesic, see~\cite{EM1970,Bruveris2014,Bruveris2017}.

To show the existence of minimizing geodesics we  use \cite[Remark~5.4]{Bruveris2014b_preprint}. There the existence of minimizing geodesics for metrics on $\mathcal I^2(S^1,\R^d)$ is proven provided that: 
\begin{itemize}
    \item they are uniformly bounded and uniformly coercive with respect to the background $\|\cdot\|_{H^n(d\theta)}$-norm on metric balls in the geodesic distance;
    \item they are of the form
\[
G_c(h,h) = \sum_{k=1}^N \| A_k(c) h \|_{F_k}^2\,,
\]
with some Hilbert spaces $F_k$ and smooth maps $A_k : \mathcal I^2 \to L(H^2, F_k)$, where the maps $A_k$ are required to have the property:
\[
c^j \to c \text{ weakly in } H^1_t \mathcal I^2_\theta
\quad\Rightarrow\quad
A_k(c^j) \dot c^j \to A_k(c) \dot c
\text{ weakly in } L^2(I,F_k)\,.
\]
\end{itemize}
In our case $H^1_t \mathcal I^2_\theta = H^1([0,1], \mathcal I^2(S^1,\R^d))$. In our case $N=4$, $F_k = L^2(S^1,\R^d)$ and $A_1(c)h = \sqrt{a_0} h$, $A_2(c)h = \sqrt{a_1} D_s h^\top$, $A_3(c)h = \sqrt{b_1} D_s h^\perp$ and $A_4(c)h = \sqrt{a_2} D_s^2 h$. The necessary convergence properties follow from \cite[Lemma~5.9]{Bruveris2014b_preprint}.
\end{proof}

For applications in matching the central task is to obtain  stable and fast algorithms to calculate the induced geodesic distance. Our framework is based on discretizing the Riemannian energy and minimizing it over all (discrete) paths. For some first order metrics there exist transformations that can significantly speed up these calculations, because they yield explicit formulas for the geodesic distance on open parametrized curves. In related work on open planar curves Kurtek and Needham \cite{KuNe2018} will follow this approach to obtain fast numerical algorithms for first order metrics. 
 The aim of the present article is to develop a numerical framework for a wider class of metrics on open and closed curves, that should allow one to model a variety of different matching behaviors. Furthermore, we plan to further enhance our framework in future work in order to be able to deal with surfaces in addition to curves. For these reasons we do not take advantage of these explicit formulas.

\subsection{Elastic metrics on shape spaces of curves} 
We will now use the metrics  defined in Section~\ref{sec:elastic_metrics_param} to induce Riemannian metrics on shape spaces of unparametrized curves, as defined in Section~\ref{sec:shape_space}. In this section we fix a subgroup $H$ of the group $SE(d)$ and let
$\mathcal S_f(M^1,\R^d)=\left(\Imm_f(M^1,\R^d)/\on{Diff}(M^1)\right)/H$. Using invariance properties of the metrics, we obtain the following result
concerning the induced metrics on the quotient space:
\begin{theorem}\label{thm:metric_shape_space} 
Let $G$  be an elastic metric of type \eqref{eq:sobolev_metric:constant_coeff} 
or \eqref{eq:sobolev_metric:scale_inv}. Then $G$ induces a Riemannian metric on the quotient space
$\mathcal S_f(M^1,\R^d)$ such that the projection 
$$\pi : \Imm_f(M^1,\R^d) \to \mathcal S_f(M^1,\R^d)$$ 
is a Riemannian submersion.  For the metric $G^2$ the result holds also for quotient spaces with respect to scale.
\end{theorem}
\begin{remark}
\normalfont 
Note, that the invariance of the metric is only a necessary condition for a metric on the space of parametrized curves to induce a metric on unparametrized curves. To guarantee the existence of the induced quotient metric one has to verify the existence of the horizontal bundle for each specific metric. This has been achieved for all the metrics studied in this article \cite{Michor2007}.
\end{remark}

We have to restrict ourselves to free immersions to obtain a smooth structure on the quotient $\mathcal S_f(M^1,\R^d)$. However, in numerical calculations we will work with the full space $\on{Imm}(M^1,\R^d)$ and the quotient 
$\mathcal S(M^1,\R^d) = \Imm(M^1,\R^d)/\on{Diff}(M^1)/H$.
This quotient space has singularities at non-free immersions, but the singularities are very mild~\cite{Michor1991}.

The geodesic distance on $\on{Imm}(M^1,\R^d)$ gives rise to a distance on the quotient space $\mathcal S(M^1,\R^d)$, which coincides---for shapes sufficiently close to each other---with the geodesic distance of the induced Riemannian metric on $\mathcal S_f(M^1,\R^d)$. Here we use the fact that $\mathcal S_f(M^1,\R^d)$ is an open dense subset of $\mathcal S(M^1,\R^d)$ and thus the geodesic distances
coincide at least as long as a minimizing deformation  does not encounter points in $\mathcal S(M^1,\R^d)$ that do not
belong to $\mathcal S_f(M^1,\R^d)$.

This distance can be then calculated using paths in $\Imm(M^1,\R^d)$ connecting $c_0$ to the orbit $c_1 \circ \on{Diff}(M^1)\circ H$, i.e., 
for $\pi(c_0),\pi(c_1) \in \Imm(M^1,\R^d)/\on{Diff}(M^1)/H$ we have,
\begin{equation}
\label{eq:dist_shape_direct}
\on{dist}\big(\pi(c_0), \pi(c_1)\big) = \inf \left\{ L(c) \,:\, c \in \mathcal P,\, c(0) = c_0,\, c(1) \in c_1 \circ \on{Diff}(M^1)\circ H\right\} \,.	
\end{equation}
We have the following completeness result for the quotient space of closed curves modulo reparametrizations:
\begin{theorem}
Let $G$ be an elastic metric of type \eqref{eq:sobolev_metric:constant_coeff} 
or \eqref{eq:sobolev_metric:scale_inv} with $a_0,a_2>0$ on the shape space
$\mathcal S(S^1,\R^d)$. Then the metric completion of $\mathcal S(S^1,\R^d)$ equipped with the quotient distance is the space of all Sobolev shapes of class $H^2$,
$$\mathcal I^2(S^1,\R^d)/\mathcal D^2(S^1)/ H\;.$$
Here $\mathcal D^2$ denotes the $H^2$-Sobolev completion of the diffeomorphism group $\on{Diff}(S^1)$.
   Furthermore the metric completion is a length space and any two shapes in the same connected component can be joined by a minimizing geodesic.
\end{theorem}

The proof of this theorem is verbatim the same as in \cite[Section~6]{Bruveris2014b_preprint}.

\section{Oriented varifold metrics}
\label{sec:VarifoldDistance}
We will derive an efficient relaxation term for the matching constraint using distances on the space of curves that originate from geometric measure theory and which have been applied extensively in shape analysis and computational anatomy. Heuristically, the philosophy is to induce a metric on the space of unparametrized curves using their representations as generalized distributions.


Several models for such distributions and their associated metrics have been proposed: measures \cite{Glaunes2004}, currents \cite{Glaunes2006,Durrleman2010,Benn2017} or varifolds \cite{Charon2013}. The most recent work \cite{Charon2017} introduces the general representation of a curve as an oriented varifold which combines the different approaches into a single framework.
Varifolds can be used, in principle, not just to define distances between curves but between embedded submanifolds of any dimension and codimension although numerical implementations exist only for curves and surfaces in $\R^2$ and $\R^3$.
In the present work we will focus only on smooth, open or closed curves. However, we have to be careful when applying the varifold framework because, technically, we are dealing with immersed and not embedded curves. In what follows, we give a brief reminder of the framework of oriented varifolds as it applies to curves while addressing the distinction between immersed and embedded curves in more detail.

\subsection{Representation of curves as oriented varifolds}
Let $C_0(\R^d \times S^{d-1})$ be the space of continuous functions vanishing at infinity. An oriented varifold is intuitively a joint distribution of point positions and directions. Mathematically, we define them as follows:
\begin{definition}
An oriented varifold is an element of the distribution space $C_0(\R^d \times S^{d-1})^*$, i.e., a signed measure on the product $\R^d \times S^{d-1}$.
\end{definition}

The analogue of Dirac distributions in the context of oriented varifolds are the distributions $\delta_{(x,u)}$ with $(x,u) \in \R^{d} \times S^{d-1}$, defined by $\delta_{(x,u)}(\omega) = \omega(x,u)$ for all test functions $\omega \in C_0(\R^d \times S^{d-1})$.

Next we define the natural representation of curves as oriented varifolds.
\begin{definition}
The varifold application $\mu: c \mapsto \mu_c$ associates to any immersion $c \in \on{Imm}(M^1,\R^d)$ the oriented varifold $\mu_c$ defined, for any $\omega \in C_0(\R^d \times S^{d-1})$, by
\begin{equation}
 \label{def_varifold}
 \mu_c(\omega) = \int_{M^1} \omega\left(c(\theta),\frac{c'(\theta)}{|c'(\theta)|}\right) \!\ud s\,.
\end{equation}
\end{definition}
Note that \eqref{def_varifold} writes informally as $\mu_c = \int_{M^1} \delta_{(c(\theta),u(\theta))} \ud s$ where $u(\theta) = \frac{c'(\theta)}{|c'(\theta)|}$ is the unit tangent vector at $\theta$; in other words $\mu_c$ can be interpreted as the weighed combination of Diracs at the point positions of the curve $c(\theta)$ with attached vectors $u(\theta)$ and infinitesimal weights given by the arclength $\ud s = |c'(\theta)| \ud \theta$. 

A key property is that $\mu_c$ is actually \emph{independent} of the parametrization. Indeed, a straightforward change of variables in \eqref{def_varifold} shows that for any positive reparametrization $\ph$ in $\text{Diff}^+(M^1)=\{\ph \in \text{Diff}(M^1) \ | \ \ph'(\theta)>0 \}$, one has $\mu_{c \circ \ph} = \mu_c$. It follows that the map $c \mapsto \mu_c$ projects to a well-defined map from the quotient space $B_{i}^+(M^1,\R^d) = \on{Imm}(M^1,\R^d)/\text{Diff}^+(M^1)$ of \emph{oriented} unparametrized immersed curves into the space of varifolds:
\[
   \begin{tikzcd}[row sep=huge,column sep=huge]
     \text{Imm}(M^1,\R^d) \arrow{r}{\pi^+} \arrow[swap]{dr}{\mu} & B_{i}^+(M^1,\R^d) \arrow{d}{} \\
      & C_0(\R^d \times S^{d-1})^*
   \end{tikzcd}
\]
Note that neither the original map $\mu_c$ nor the resulting quotient map are surjective as there are many varifolds that are not curves (e.g. a single Dirac). Whether the quotient map is injective is a question with a more nuanced answer and we will discuss it more thoroughly in the following.
\begin{remark}
\normalfont
 The varifold representation remains sensitive to the orientation of a curve because, if $\check{c}$ is the curve $c$ with the opposite orientation then we have, in general, $\mu_{\check{c}} \neq \mu_{c}$. In all of the applications considered in this paper curves can be naturally and consistently oriented and hence orientation represents a relevant piece of information that can be exploited for curve matching (cf the discussion in \cite{Charon2017}). However, as we will explain in the next subsection, we can also consider the quotient spaces of unoriented curves and varifolds by constraining test functions to be symmetric with respect to the second variable $u$; this corresponds to the framework of unoriented varifolds of \cite{Charon2013}.       
\end{remark}

\subsection{Kernel metrics}
\label{ssec:varifold_kernel} 
The varifold application embeds unparametrized curves in a common space of distributions. This suggests that we can construct distances on the space of curves by restricting distances or pseudo-distances defined on the space of varifolds. The natural choice would be the distance induced by the norm on $C_0(\R^d, S^{d-1})^\ast$ that is dual to the supremum norm on $C_0(\R^d \times S^{d-1})$. However, this yields a fundamentally nonsmooth distance for curves, as it essentially measures the exact overlap between two curves. Thus, to obtain more reasonable distances, one needs to restrict oneself to more regular spaces of test functions equipped with stronger norms. 

A practical approach is to consider a Hilbert space $\mathcal H$ of test functions, continuously embedded in $C_0(\R^d, S^{d-1})$. In this case $\mathcal H$ is a \emph{reproducing kernel Hilbert space (RKHS)} and it is generated by a positive definite kernel $k$ on the product space $\R^d \x S^{d-1}$. Following \cite{Charon2017}, we require the kernel $k$ to have some additional structure, namely, $k$ has to be a product of a radial kernel $\rho$ on $\R^d$ and a zonal kernel $\gamma$ on $S^{d-1}$, i.e.,
\begin{equation}
\label{eq:sep_kernel}
k(x,u,y,v) \doteq \rho(|x-y|^2)\gamma(u \cdot v)
\end{equation}
for all $(x,u)$ and $(y,v)$ in $\R^d \times  S^{d-1}$. Here $\rho$ defines a continuous, positive radial basis function with $\rho(t) \rightarrow 0$ as $|t| \rightarrow \infty$ and $\gamma : [-1,1] \to \R$ defines a continuous zonal function on the sphere. The general theory of reproducing kernels \cite{Aronszajn1950} states that under these assumptions on the kernel $k$, the Hilbert space $\mathcal H$ is uniquely determined by $k$. We denote by $\langle \cdot,\cdot \rangle_{\mathcal H}$ the inner product on $\mathcal H$.
The Riesz duality map then induces an inner product on $\mathcal{H}^*$. Moreover, using the dual map $C_0(\R^d \times S^{d-1})^* \rightarrow \mathcal{H}^*$, it then also induces a pseudo-distance on the space of varifolds (which is actually a distance if $\mathcal H$ is additionally assumed to be dense in $C_0(\R^d \times S^{d-1})$). Having chosen a kernel of the form \eqref{eq:sep_kernel}, we will use the generic notation $\langle \cdot, \cdot \rangle_{\on{Var}}$ for the associated inner product $\langle\cdot,\cdot\rangle_{\mathcal H^*}$ on the space of varifolds. 

For the purpose of this paper, we are in fact interested in the metric that is induced on the space of curves by the varifold map $c \mapsto \mu_c$. This metric is given by $d^{\on{Var}}(c_1,c_2) = \|\mu_{c_1} - \mu_{c_2}\|_{\on{Var}} = \langle \mu_{c_1} - \mu_{c_2}, \mu_{c_1} - \mu_{c_2} \rangle^{1/2}_{\on{Var}}$.  The reproducing kernel property implies---cf. \cite{Charon2017} for details---that for any two curves $c_1,c_2$ we have 
\begin{equation}
 \label{eq:metric_W_curves}
 \langle \mu_{c_1} , \mu_{c_2} \rangle_{\on{Var}} = \iint_{M^1\x M^1} \rho(|c_1(\theta_1) - c_2(\theta_2)|^2) \gamma\left(\frac{c_1'(\theta_1)}{|c_1'(\theta_1)|} \cdot \frac{c_2'(\theta_2)}{|c_2'(\theta_2)|} \right) \ud s_1 \ud s_2\,,
\end{equation}
leading to a similar closed-form expression for $d^{\on{Var}}(c_1,c_2)^2$. Equation \eqref{eq:metric_W_curves} shows that $d^{\on{Var}}$ can be interpreted as a localized comparison between the relative positions of points and tangent lines of the two curves, quantified by the choice of kernel functions $\rho$ and $\gamma$. We will see later how to efficiently evaluate  these expressions numerically for discrete curves. 

In general $d^{\on{Var}}$ only defines a pseudo-distance. In order to be able to separate any two curves, one needs to ensure that the space of test functions $\mathcal{H}$ is large enough. A sufficient condition for this is given by the following theorem which is a particular case of Proposition 4 in \cite{Charon2017}:
\begin{theorem}
\label{thm:vardistance}
Assume that $\rho$ and $\gamma$ are $C^1$-functions, $\rho$ is $C_0$-universal and $\gamma(1)>0$. Then, if $d^{\on{Var}}(c_1,c_2)=0$, we have $\on{Im}(c_1)=\on{Im}(c_2)$ i.e., the images of $c_1$ and $c_2$ in $\R^d$ coincide. 
\end{theorem}
The kernel $\rho$ is said to be \emph{$C_0$-universal} if its associated reproducing kernel Hilbert space is dense in $C_0(\R^d,\R)$. This is the case for Gaussian, Cauchy and Wendland kernels for example, we refer to \cite{Carmeli2010} or \cite{Sriperumbudur10} for details on the construction and characterization of such kernels.

Note that while we have $\text{Im}(c_1)=\text{Im}(c_2)$ in the result of Theorem \ref{thm:vardistance}, it is not necessarily the case that the orientations of $c_1$ and $c_2$ coincide. This is in particular not true when $\gamma$ defines an orientation-invariant kernel (cf discussion below). This can be enforced, however, under the following conditions: 
\begin{corollary}
\label{cor:vardistance_embeddings}
In addition to the assumptions of Theorem \ref{thm:vardistance}, if the function $\gamma$ is such that $\gamma(-t) \neq \gamma(t)$ for all $t \neq 0$, then $d^{\on{Var}}$ defines a distance on the space of oriented, unparametrized, embedded curves.  
\end{corollary}
The proof follows from very similar arguments than the one of Theorem \ref{thm:vardistance} in \cite{Charon2017} that we do not repeat for concision. Note that the last statement is equivalent to saying that the varifold application $\mu$ into $\mathcal{H}^*$ is injective if restricted to $\on{Emb}(M^1,\R^d)/\text{Diff}^+(M^1)$. 


\begin{figure}
    \centering
    \begin{tabular}{ccc}
    \includegraphics[width=5cm]{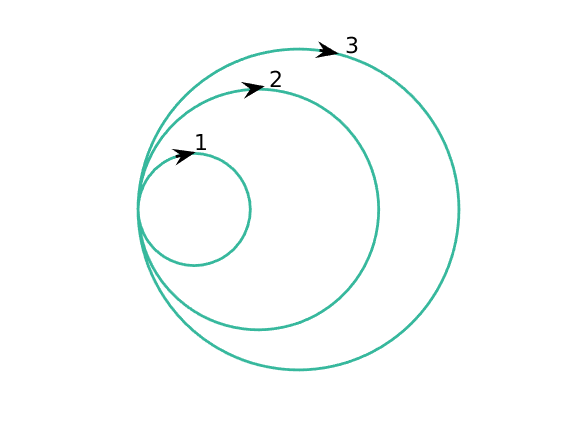} & & \includegraphics[width=5cm]{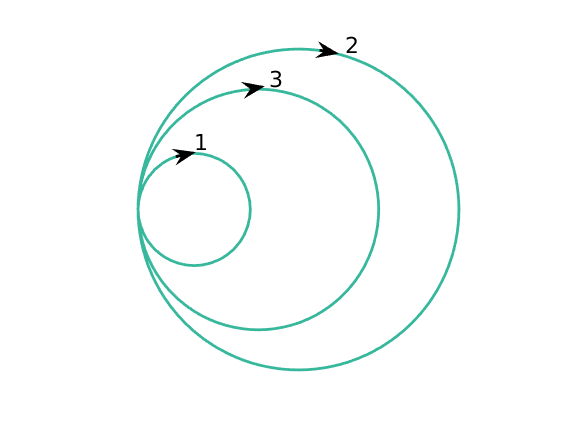}
    \end{tabular}
    \caption{Example of two distinct oriented immersed curves that are however equal in the space of oriented varifolds. The numbers reflect the order of crossing of each piece.}
    \label{fig:three_circles}
\end{figure}

Unfortunately, it is easy to see that this result cannot be extended to the larger space of free immersions. Indeed, the varifold representation $\mu_c$ of a curve $c$ takes into account only the image of an immersed curve and its orientation. The example in Figure~\ref{fig:three_circles} (also mentioned in \cite{Benn2017} for the specific case of currents) shows two distinct elements of $B_{i}^+(S^1,\R^d) = \on{Imm}(S^1,\R^d)/\text{Diff}^+(S^1)$, both projections of free immersions, that are nevertheless equal as oriented varifolds. Yet, the previous result on embedded curves can be generalized to the space of immersed curves with a finite number of transverse self-intersections. Note that any such immersion is already a free immersion.

\begin{theorem}
\label{thm:vardistance_immersions}
With the assumptions of Corollary \ref{cor:vardistance_embeddings}, if $c_1$ and $c_2$ are two immersions with a finite number of transverse self-intersections then $d^{\on{Var}}(c_1,c_2)=0$ if and only if the two curves coincide in $B_{i}^+(M^1,\R^d)$.
\end{theorem}
\begin{proof}
 Let $c_1$ and $c_2$ be two such immersions with $d^{\on{Var}}(c_1,c_2)=0$. We need to show that there exists $\ph \in \text{Diff}(M^1)$ such that $c_2=c_1 \circ \ph$. 
 
 Using Theorem \ref{thm:vardistance} we already know that $\text{Im}(c_1)=\text{Im}(c_2)$. Let us denote by $(p_1,\ldots,p_m) \in (\R^d)^m$ the self-intersection points of $c_1$, which are the same for $c_2$ as the two curves have the same image. For any $x \in \text{Im}(c_1) \backslash \{p_1,\ldots,p_m\}$, there exist unique preimages $c_1^{-1}(x)$ and $c_2^{-1}(x)$ in $M^1$. Let us denote by $\Theta_1^s$ and $\Theta_2^s$ the preimages of the self-intersection points under $c_1$ and $c_2$ respectively. We set $\ph(\theta)=c_1^{-1} \circ c_2(\theta)$ for $\theta \in M^1 \backslash \Theta_2^s$. 
 
 First, note that $\ph$ is smooth on $M^1 \backslash \Theta_2^s$ and that $c_2(\theta) = c_1 \circ \ph(\theta)$ for all $\theta \in M^1 \backslash \Theta_2^s$. In addition, using Corollary \ref{cor:vardistance_embeddings}, we also have that the orientation coincide on each connected component of $c_1(M^1 \backslash \Theta_2^s)$ and therefore $\ph'(\theta)>0$ on $M^1 \backslash \Theta_2^s$. 
 
 Now, let $i \in \{1,\ldots,m\}$ and $c_2^{-1}(p_i) = \{\theta_1^i,\ldots,\theta_{m_i}^i\}$. Since all self-intersections are transverse, we have that the $c_2'(\theta_k^i)$ are all distinct from one another. Similarly, as the two curves have the same image, we can write  $c_1^{-1}(p_i) = \{\tilde{\theta}_1^i,\ldots,\tilde{\theta}_{m_i}^i\}$ and with the adequate ordering we can also assume that the two vectors $c_1'(\tilde{\theta}_k^i)$ and $c_2'(\theta_k^i)$ are collinear. Then, setting for all $i$ and $k$, $\ph(\theta_k^i) = \tilde{\theta}_k^i$, we obtain a bijection $\ph: \ M^1 \rightarrow M^1$ such that $c_2 = c_1 \circ \ph$. Moreover, the above construction makes $\ph$ a smooth function that satisfies $c_2'(\theta) = c_1'(\ph(\theta)) \ph'(\theta)$ for all $\theta \in M^1$.  
\end{proof}

In certain situations, it may be more relevant to work with unoriented immersed curves i.e. with the space $B_{i}(M^1,\R^d) = \on{Imm}(M^1,\R^d) / \on{Diff}(M^1)$. An equivalent result holds by instead considering orientation-invariant kernels which are such that $\gamma(-t)=\gamma(t)$ for all $t\in [-1,1]$. Indeed one can easily see from \eqref{eq:metric_W_curves} that the resulting metric $d^{\on{Var}}$ is also invariant to orientation changes in either curve; this corresponds to the subclass of unoriented varifold metrics, c.f. \cite{Charon2013,Charon2017}. Theorem \ref{thm:vardistance_immersions} can be then replaced by:
\begin{corollary}
\label{cor:vardistance_immersions}
With the assumptions of Theorem \ref{thm:vardistance} and if $\gamma$ is an even function, two immersed curves with finite numbers of transverse self-intersections are equal in $B_{i}(M^1,\R^d) = \on{Imm}(M^1,\R^d) / \on{Diff}(M^1)$ if and only if we have $d^{\on{Var}}(c_1,c_2)=0$. 
\end{corollary}

In summary, although the varifold metrics introduced here may not always distinguish two given immersed curves, we see from Theorem \ref{thm:vardistance_immersions} and Corollary \ref{cor:vardistance_immersions} that this will only occur in pathological situations such as shown in Figure \ref{fig:three_circles}. We will typically ignore such cases in the practical curve matching applications of this paper.

\subsection{Varifold distance as a constraint}
 \label{sec:InexactMatching}
The invariance of varifold-induced distances under reparametrizations makes them a natural tool for enforcing the exact matching constraint in the geodesic boundary value problem for elastic metrics. Indeed, the geodesic distance $\on{dist}(\pi(c_0), \pi(c_1))^2$ can be computed  in the following way,
\begin{equation}\label{Energy_Vari}
\on{dist}(\pi(c_0), \pi(c_1))^2 = 
\inf\left\{ E(c) \,:\, c \in \mathcal P,\, c(0) = c_0,\, d^{\on{Var}}(c(1),c_1)^2=0 \right\}\,,
\end{equation}
where $E(c)$ is the Riemannian energy of the path $c$ and $\mathcal P$ is the space of all smooth paths in $\on{Imm}(M^1,\R^d)$. The squared varifold distance, which can be calculated explicitly via \eqref{eq:metric_W_curves}, is used as a smooth constraint enforcing the endpoint condition $\pi(c(1)) = \pi(c_1)$. In contrast with the direct approach of calculating $\on{dist}(\pi(c_0),\pi(c_1))^2$ via \eqref{eq:dist_shape_direct}, the formulation \eqref{Energy_Vari} does not require optimization over reparametrizations. Never the less the optimal point correspondences can be inferred from our method. The equivalence between the two formulations is rigorous provided the curves $c_0$, $c_1$ and the kernel $k$ satisfy the assumptions of Corollary \ref{cor:vardistance_immersions}. In that case we have $d^{\on{Var}}(c(1),c_1)^2 = 0 \Leftrightarrow c(1)=c_1 \circ \ph$ for some $\ph \in \text{Diff}(M^1)$. 
\begin{remark}
\normalfont
Note that this corresponds to the problem of matching unoriented, unparametrized, immersed curves, i.e., elements of $B_{i}(M^1,\R^d)$. In certain other situations, one could assume that curves have been consistently oriented from the start and wish to solve the matching problem for oriented curves in $B_{i}^+(M^1,\R^d)$ instead. In that case, it is not difficult to see that we can also reformulate the problem as \eqref{Energy_Vari} by choosing an orientation-sensitive metric for $d^{\on{Var}}$.   
\end{remark}


\subsection{Invariance to similarities}
As mentioned at the end of Section \ref{sec:shape_space}, it is often important to compare curves modulo the positive similarity group $S(d)$ and therefore quotient out these transformations in the estimation of distance and geodesic. We first focus on the particular subgroup $SE(d)$ of Euclidean motions. In that case, we have seen that both families of elastic metrics $G_c^1$ in \eqref{eq:sobolev_metric:constant_coeff} and $G_c^2$ in \eqref{eq:sobolev_metric:scale_inv} are invariant to the action of $SE(d)$. Thanks to the particular form of $k$ in \eqref{eq:sep_kernel}, it turns out that this is also the case of the kernel-based distances $d^{\on{Var}}$, i.e, we always have $d^{\on{Var}}(A.(c_1+w),A.(c_2+w)) = d^{\on{Var}}(c_1,c_2)$ for any two curves $c_1,c_2$ and rigid motion $(A,w) \in SE(d)$. Consequently, the invariant matching problem: 
\begin{equation*}
\inf \left\{ E(c):\, c \in \mathcal P,\, (A,w) \in SE(d), \, c(0) = c_0,\, c(1)=A.(c_1+w)\right\}
\end{equation*}
becomes once again equivalent to
\begin{equation*}
\inf \left\{ E(c) \,:\, c \in \mathcal P,\, (A,w) \in SE(d), \, c(0) = c_0,\, d^{\on{Var}}(c(1),A.(c_1+w))^2=0\right\}.	
\end{equation*}
which we can then tackle like previously in either the relaxed or augmented Lagrangian formulation, jointly over the path $c$ and the finite-dimensional variables $(A,w)$. 

The case of scale-invariance is however more involved in the present setting. While the second family of elastic metrics $G_c^2$ is invariant to rescaling, this is not true for the oriented varifold metrics of Section \ref{ssec:varifold_kernel}. In fact, it is quite easy to see that no metric originating from a kernel of the form of \eqref{eq:sep_kernel} is scale-invariant as this would impose that $\rho(\lambda^2 t)= \rho(t)/\lambda^2$ for all $\lambda$ and $t$ and thus lead to a singularity at $0$ for the function $\rho$. In most applications \cite{Charon2013,Charon_thesis,Charon2017}, it is rather customary to specify kernels with an intrinsic notion of scale by setting for instance the kernel defined by $\rho$ to be a Gaussian $\rho(|x-y|^2)=e^{-\frac{|x-y|^2}{\sigma^2}}$ or a sum of Gaussian for multiscale applications. In the context of this work, we point that out the lack of invariance of $d^{\on{Var}}$ to rescaling will not constitute an issue as the varifold metric is only used as a surrogate for the matching constraint of the two immersed curves, and thus only invariance with respect to reparametrizations is necessary. 

\section{Implementation}
\label{sec:numerics}
In this section we will describe how to discretize and solve the constrained optimization problem \eqref{Energy_Vari} using both an inexact one-shot method and a iterative augmented Lagrangian scheme which enforces a better constraint satisfaction. Our code is available on GitHub\footnote{\url{https://www.github.com/h2metrics/h2metrics}}.
\subsection{A B-spline discretization}

In order to evaluate the energy functional \eqref{eq:EnergyFunctional} and the constraint \eqref{eq:metric_W_curves} we discretize paths of curves using tensor product B-splines on knot sequences of orders $n_t$ in time and $n_\theta$ in space (typically we choose $n_t=2$ and $n_\theta=3$). This produces $N_t \x N_\theta$ basis splines, with $N_t$ and $N_\theta$ being the number of control points in each variable respectively (typical values we shall take in the experimental section are $N_{\theta}=100$ and $N_t=10$), and we can write
\begin{equation}
\label{eq:TensorPath}
c(t, \th) = \sum_{i=1}^{N_t} \sum_{j=1}^{N_\th} c_{i,j} B_i(t) C_j(\th)\,.
\end{equation}
Here $B_i(t)$ are B-splines defined by an equidistant simple knot sequence on $[0,1]$ with full multiplicity at the boundary knots, and $C_j(\theta)$ are defined by an equidistant simple knot sequence on $[0,2\pi]$ with periodic boundary conditions or full multiplicity at the boundary for closed or open curves respectively; for details see Section 3 of \cite{BBHM2017}. The full multiplicity of the boundary knots in $t$ implies
\begin{align*}
c(0, \th) &= \sum_{j=1}^{N_\th} c_{1,j} C_j(\th)\,, &
c(1, \th) &= \sum_{j=1}^{N_\th} c_{N_t,j} C_j(\th)\,.
\end{align*}
Thus the initial curve $c(0)$ is given by the control points $c_{1,j}$ only, which we can utilize later for the constraint satisfaction. In terms of the standard differential operator $\partial_{\theta}$ (as opposed to arc-length differentiation $D_s$) 
the Riemannian metrics \eqref{eq:sobolev_metric:constant_coeff} and \eqref{eq:sobolev_metric:scale_inv} read as
\begin{align*}
     &G^1_c(h,k)= \int_0^{2\pi} a_0 |c'| \langle h, k \rangle + \frac{a_1}{|c'|}  \langle h'^\top, k'^\top \rangle + \frac{b_1}{|c'|}  \langle h'^\bot, k'^\bot \rangle \\
&\qquad
+ \frac{a_2}{|c'|^7} \langle c', c'' \rangle^2 \langle h', k' \rangle
- \frac{a_2}{|c'|^5} \langle c', c'' \rangle \big( \langle h', k'' \rangle + \langle h'', k' \rangle \big) +
\frac{a_2}{|c'|^3} \langle h'', k'' \rangle
\ud \th\,. \\
    &G^2_c(h,k)= \int_0^{2\pi} \frac{a_0}{\ell^3} |c'| \langle h, k \rangle + \frac{a_1}{\ell|c'|} \langle h'^\top, k'^\top \rangle + \frac{b_1}{\ell|c'|} \langle h'^\bot, k'^\bot \rangle \\
&\qquad
+ \frac{\ell a_2}{|c'|^7} \langle c', c'' \rangle^2 \langle h', k' \rangle
- \frac{\ell a_2}{|c'|^5} \langle c', c'' \rangle \big( \langle h', k'' \rangle + \langle h'', k' \rangle \big) +
\frac{\ell a_2}{|c'|^3} \langle h'', k'' \rangle
\ud \th\,.
\end{align*}
Plugging these expressions into \eqref{eq:EnergyFunctional} gives an explicit expression, which we leave out, for the energy of a given path. For a B-spline path, we approximate the integrals in the energy functional \eqref{eq:EnergyFunctional} and varifold distance \eqref{eq:metric_W_curves} using Gaussian quadrature with quadrature sites placed between knots where the curves are smooth. This yields a fast and robust way to evaluate the energy of paths. The same is true for the evaluation of the derivatives found in Appendix \ref{sec:Derivatives}. 

\subsection{The optimization procedures}
We will now describe two different methods to approximately factor out the action of the diffeomorphism group. At the end we will comment on the action of euclidean motions and scalings.

In section \ref{sec:InexactMatching} we showed that in order to factor out the diffeomorphism group we have to solve an optimization problem under the constraint that the end point of the curve satisfies $d^{\on{Var}}(c(1),c_1)^2 = 0$; this corresponds to an \textit{exact} matching of the end point and the target curve. Inspired by the paradigms of other methods like LDDMM, as a simple method we consider an \textit{inexact} matching problem where we only desire that the constraint violation of the end point is small instead of requiring it to be exactly zero. To this end a fixed large value of $\lambda$ is chosen and the following relaxed quadratic penalty functional is considered
\begin{equation}\label{eq:Energy_Vari}
\inf \left\{ E(c) +\lambda d^{\on{Var}}(c(1),c_1)^2
\,:\, c \in \mathcal P,\, c(0)=c_0\right\}\,.
\end{equation}
Here $\lambda > 0$ is a balance parameter between the elastic energy and the varifold fidelity term. This is particularly well-suited to noisy situations in which exact matching might lead to irrelevant solutions. Note that exact matching is still theoretically recovered in the limit $\lambda \rightarrow +\infty$. To solve the unconstrained optimization problem we use the HANSO library \cite{Hanso}, which utilizes a L-BFGS method. This approach does not yield a geodesic with the correct endpoint but with an appropiate choice of $\lambda$ the varifold distance term is small in practice.  In \cite{BauerMICCAI2017} we employed this method, but the problem seemed quite sensitive to the choice of $\lambda$: too small and a bad matching is achieved; too big and the optimization algorithm has difficulty finding a solution. 

In order reduce the sensitivity of the solution to the choice of the weight parameter $\lambda$ and to possibly solve the exact matching problem, we also propose an augmented Lagrangian scheme. In practice it will not be feasible for a B-spline path to satisfy $d^{\on{Var}}(c(1),c_1) = 0$ exactly, hence we would rather relax the constraint to an inequality
\begin{equation*}
    d^{\on{Var}}(c(1),c_1) \leq \varepsilon,
\end{equation*}
for some small chosen constraint error tolerance $\varepsilon > 0$. In order to solve this inequality constrained minimization  we use a simple adaptation of the augmented Lagrangian scheme, see \cite{Nocedal2006}. We introduce the augmented Lagrangian functional
\begin{equation}\label{eq:Energy_AugLag}
\mathcal{L}(c,\lambda,\mu) = E(c) - \lambda d^{\on{Var}}(c(1),c_1)^2 + \frac{\mu}{2} d^{\on{Var}}(c(1),c_1)^4\,;
\end{equation}
here $\lambda$ plays the role of the (real-valued) Lagrange multiplier associated to the constraint $d^{\on{Var}}(c(1),c_1)^2 = 0$. Notice that if $\mu = 0$ then the functional is the same as the quadratic penalty with the sign of $\lambda$ flipped. In general this method should be better conditioned than the quadratic penalty method, and convergence can be guaranteed for the penalty parameter $\mu$ above a finite threshold, and not only for $\mu \to \infty$. The constrained problem can be then solved by simultaneously minimizing $\mathcal{L}$ over $c$ while updating the Lagrange multiplier $\lambda$. We approximately solve the sequence of unconstrained problems given by
\begin{equation}
    c_k = \underset{c \in \mathcal P,\, c(0)=c_0}{\on{argmin}} \mathcal{L}(c, \lambda_k, \mu_k)\,,
\end{equation}
where $\mu_k$ is a given sequence of positive scalars which weights the constraint error penalty term, $\lambda_k$ is the current estimate of the Lagrange multiplier which is updated via the rule
\begin{equation}
    \lambda_{k+1} = \lambda_k - \mu_k d^{\on{Var}}(c_k(1),c_1)^2\,.
\end{equation}
 At each iteration we check if the soft constraint is satisfied, if so we accept the current value of $\mu$ and continue, if not we increase $\mu$ in order to enforce the constraints. At each iteration step, we need to solve an unconstrained minimization problem, for this we use the L-BFGS method in the HANSO library \cite{Hanso}. In practice we only need to solve the sequence of problems with a sequence of gradient tolerances $\tau_k \to 0$. For small $k$ the tolerance can be chosen quite high to quickly terminate the optimization algorithm. The whole method is summarized in Algorithm \ref{alg:AugLagrangian}. As opposed to the quadratic penalty method, the augmented Lagrangian method seemed less sensitive to the choice of sequence of $\mu_k$ but at the cost of solving several unconstrained optimization problems instead of a single one. If solving each unconstrained optimization problem is difficult, it might  be computationally inefficient to use an augmented Lagrangian method.
\begin{algorithm}
\caption{Augmented Lagrangian} \label{alg:AugLagrangian}
\begin{algorithmic}
\STATE \textbf{Input:} Curves $c_0,c_1$ to be matched.
\STATE{Set $\mu_0 >0, \lambda_0 \leq 0, \tau_0 >0, \tau_{final} > 0, c_{init}^{0}(t,\theta) = c_0(\theta)$}
\FOR{$k=0,1, 2, \dots ,k_{\max}$}
\STATE $c_k = \underset{c}{\on{argmin}} \, \mathcal{L}(c,\lambda_k,\mu_k)$, with stopping criteria $\| \nabla_c \mathcal{L}(c_k,\lambda_k,\mu_k) \| < \tau_k $.
\IF{ $d^{\on{Var}}(c_k(1),c_1)^2 \leq \varepsilon$ and $\tau_k \leq \tau_{final} $ }
\RETURN $c_k$
\ENDIF
\STATE $\lambda_{k+1} = \lambda_k - \mu_k d^{\on{Var}}(c_k(1),c_1)^2$
\IF{$\| d^{\on{Var}}(c_k(1),c_1)^2 \| < \varepsilon$}
\STATE $\mu_{k+1} = \mu_k$ 
\ELSE
\STATE $\mu_{k+1} = \varrho \mu_k$ 
\ENDIF
\IF{ $\tau_k < \tau_{final}$}
\STATE $\tau_{k+1} = \frac{1}{2} \tau_k$
\ELSE
\STATE $\tau_{k+1} = \tau_k$  
\ENDIF
\STATE $c_{init}^{k+1} = c_k$ 
\ENDFOR
\end{algorithmic}
\end{algorithm}

In order to additionally factor out the action of the Euclidean motion group $\SE(d)$ and scalings, 
we can simply replace the constraint terms involving the varifold distance by
\begin{equation}
\label{eq:var_rigid_invariant}
    d^{\on{Var}}(c(1),rA(c_1+b)), \quad (r,A,b) \in \R^+ \ltimes \SE(d)
\end{equation}
and add $(r,A,b)$ to the list of variables in each unconstrained minimization subproblem. Observe that there are several orderings of the group actions that would have been possible, we choose to translate first in order to be able to center the curves before rotating them. Finally we want to add some remarks on alternative methods:

\begin{remark}[Discretization of $\on{Diff}(M)$]
\normalfont
To solve the geodesic boundary value problem on shape space, we have proposed in \cite{BBHM2017} a method that also discretizes the reparametrization group $\on{Diff}(S^1)$ using B-splines. The action of the reparametrization group is by composition, which does not preserve the B-spline space, as degrees are added. To overcome this we added an $L^2$ projection step, the composition $c \circ \ph$ is projected back into a fixed lower order spline space. This has the disadvantage that the projection can smooth out details of the original curve, depending on how many control points are used and which parts of the curve are reparametrized. Furthermore, this methods requires a good choice of an initial path, which turned out to be a nontrivial obstacle in examples where the shapes under consideration are sufficiently different from each other. The inexact  matching algorithm presented in this paper does not have these problems  as we can always choose to initialize the optimization procedure with the constant path..
\end{remark}
\begin{remark}[Dynamic Programming]
\normalfont
In the SRVF framework \cite{Jermyn2011,Esl2014b}, a dynamic programming method is usually used to find a global solution of the geodesic boundary value problem on shape space. This relies heavily on access to a fast formula for the geodesic distance between parametrized curves (for the SRVF metric and open curves there even exists an explicit analytic formula). For the bigger class of metrics considered in this article calculating the geodesic distance between parametrized curves is
comparable fast to calculating the geodesic distance between unparametrized curves. 
Thus dynamic programming, which relies on iteratively calculating geodesic distances between parametrized curves,
is not well-suited for this particular class of metrics. 
\end{remark}


\section{Experiments}

\begin{figure}
\centering
\includegraphics[width=\textwidth]{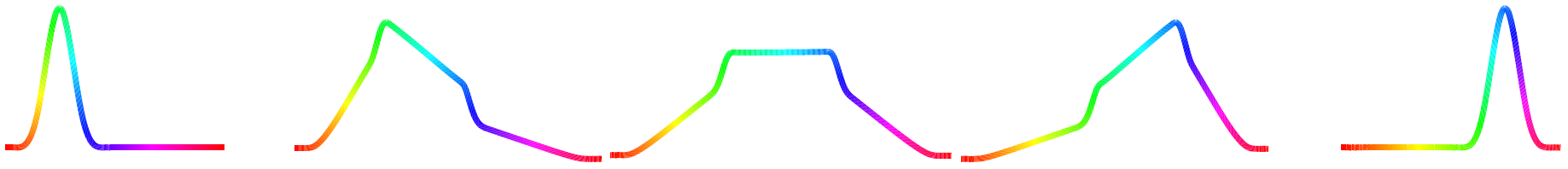} \\
\includegraphics[width=\textwidth]{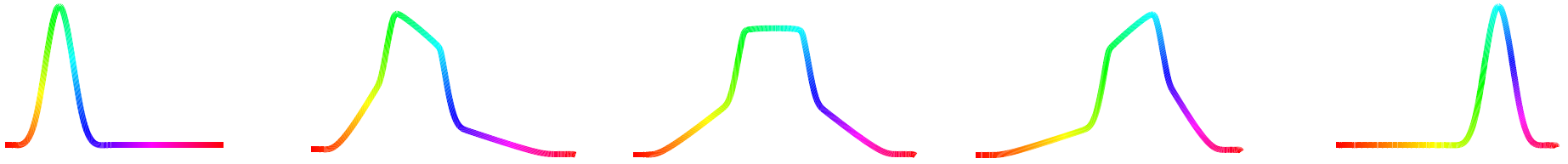} \\
\includegraphics[width=\textwidth]{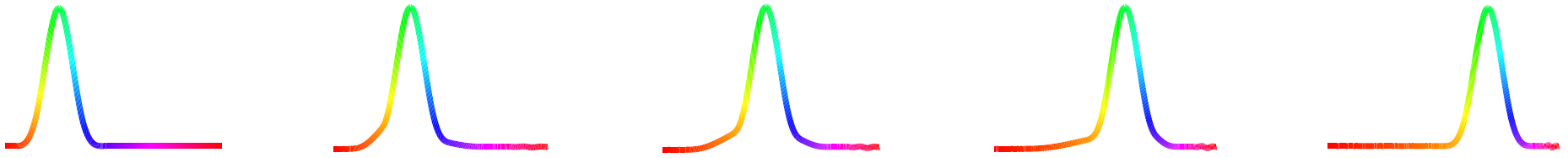}
\caption{Matching of curves with different constants in the metric. The initial and final curves are the same in all rows. The constants are
$a_0=0$, $a_1=10$ and $a_2=10^{-3}$ for all rows. The constant $b_1$ is
$0.1$ in row one, $1$ in row two and $10$ in row 3. 
} \label{fig:bumps}
\end{figure}

The choice of constants in the metric does matter---it influences both the path of the minimal geodesic as well as the parametrization of the endpoint. To illustrate this we computed geodesics between the two curves $c_0, c_1$ shown in the first and last column of Figure~\ref{fig:bumps}. Both curves have length $2\pi$ and are parametrized by arc length. We use the metric
\[
G_c(h,k) = 
\int_I a_1 \langle D_s h^\top,D_s k^\top \rangle 
+b_1\langle D_s h^{\bot},D_s k^{\bot} \rangle 
+a_2\langle D_s^2 h,D_s^2 h \rangle \ud s \;,
\]
with the following choices of constants:
\begin{align*}
(a_1, b_1, a_2) &= (10, 0.1, 10^{-3}) \\
(a_1, b_1, a_2) &= (10, 1, 10^{-3}) \\
(a_1, b_1, a_2) &= (10, 10, 10^{-3}) \,.
\end{align*}
In other words we change the relative weighting of the normal and tangential components in the $H^1$-part of the metric. The geodesics are computed modulo translations. 

We see that the geodesic in the first row bends the curve to flatten the bump in the initial curve, $c_0$, and to create the bump in the target curve, $c_1$. Note that the tip of the bump in $c_0$ is being matched to the fold at the bottom of $c_1$.  However, in the third row the bump is translated from $c_0$ to $c_1$ resulting in stretching and compression along the curve; in particular the tip of $c_0$ is matched to the tip of $c_1$. This is expected, because in this example $b_1 = 10$ and thus bending is more costly now. The middle row show intermediate behavior between the both extremes.

\subsection{Scale-invariant metrics} 
In  applications to shape analysis the scale of the curves often has  no natural meaning. To factor out scale differences, it has been proposed to re-scale the curves to fixed length and to perform the analysis on these constant length curves, see \cite{Jermyn2011,BBHM2017}. However, for a non scale-invariant metric the choice of scale, i.e., the common length of the curves, has a large effect on the resulting analysis, as demonstrated in Table~\ref{tab:scale}.

\begin{table}
\setlength\extrarowheight{3pt}
\centering
\begin{tabular}{@{}lcccccc@{}}
\toprule
& \multicolumn{2}{ c }{Fish}&\multicolumn{2}{ c }{Corp. Cal.}&\multicolumn{2}{ c }{Corp. Cal.--Fish}\\
\cmidrule{2-3} \cmidrule{4-5} \cmidrule{6-7}
&$G^1$ &$G^2$ & $G^1$& $G^2$&$G^1$ &\raisebox{-1pt}{$G^2$}\\
\midrule
$d(c_1,c_2)$&0.2423          &0.2411  &0.1305           &0.1305 &0.4367   &0.4437        \\
$d(2c_1,2c_2)$&  0.2875& 0.2411   & 0.1551 &0.1305 &  0.5231     &0.4437 \\
$d(3c_1,3c_2)$&0.3202 & 0.2411   & 0.1670 &0.1305& 0.5538    &0.4436  \\
$d(5c_1,5 c_2)$&0.3521& 0.2398   &0.1755 &0.1305 &0.5699 & 0.4435 \\
\bottomrule
\end{tabular}
\vspace{5pt}
\caption{First 4 rows: Geodesic distance for scale and non scale invariant metrics between different shapes on varying scales. The constants in the metric were chosen to be $a_0=a_1=b_1=1$,
$a_2=10^{-4}$. 
Note that one has to adapt the parameters for the varifold distance to the curve lengths. In the above table Fish refers to shapes from the Surrey fish dataset and Corp. Cal. refers to shapes from a collection of outlines of corpi callosi.}
\label{tab:scale}
\end{table}

\begin{figure}
\centering
\setlength{\tabcolsep}{0pt}
\begin{tabular}{cccc}
\includegraphics[trim = 45mm 15mm 45mm 15mm ,clip,width=3cm]{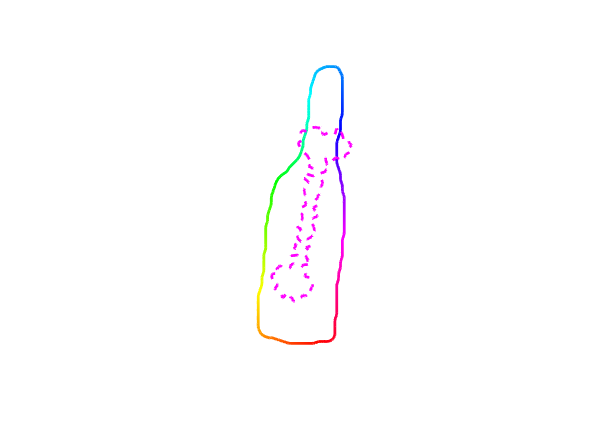} &
\includegraphics[trim = 45mm 15mm 45mm 15mm ,clip,width=3cm]{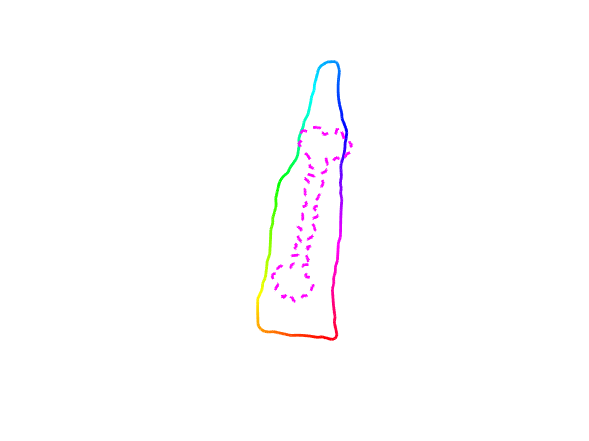} &
\includegraphics[trim = 45mm 15mm 45mm 15mm ,clip,width=3cm]{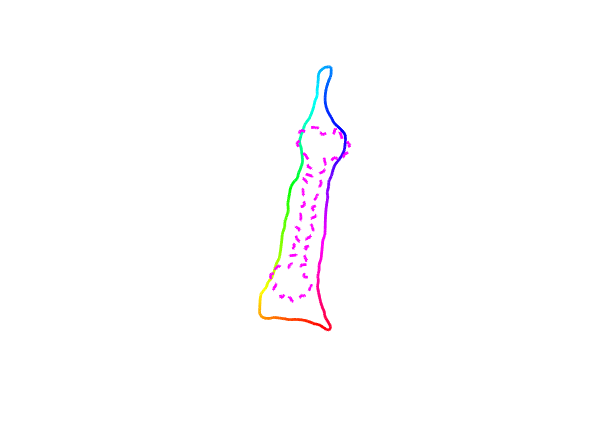} &
\includegraphics[trim = 45mm 15mm 45mm 15mm ,clip,width=3cm]{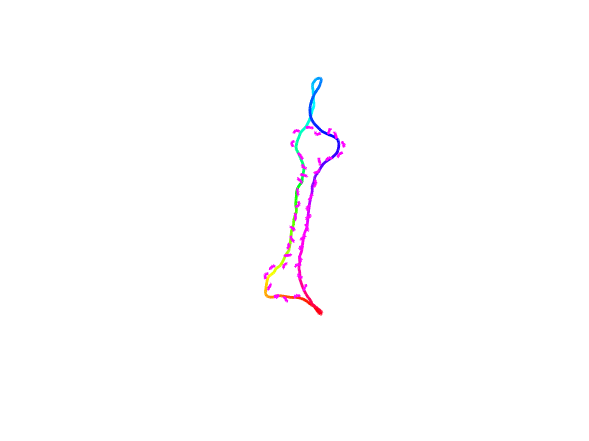}\\
t=0 & t=0.4 & t=0.7 &t=1 \\
\includegraphics[trim = 45mm 15mm 45mm 15mm ,clip,width=3cm]{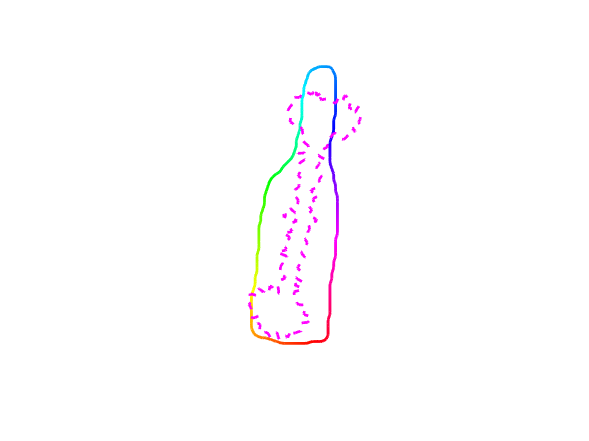} &
\includegraphics[trim = 45mm 15mm 45mm 15mm ,clip,width=3cm]{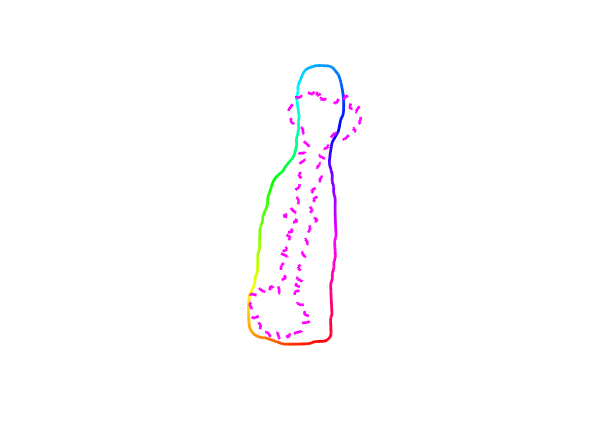} &
\includegraphics[trim = 45mm 15mm 45mm 15mm ,clip,width=3cm]{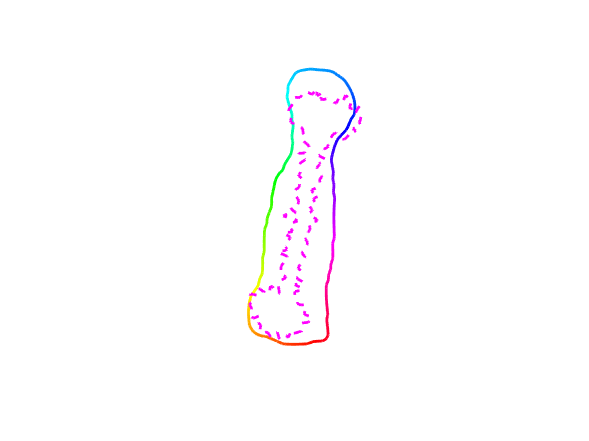} &
\includegraphics[trim = 45mm 15mm 45mm 15mm ,clip,width=3cm]{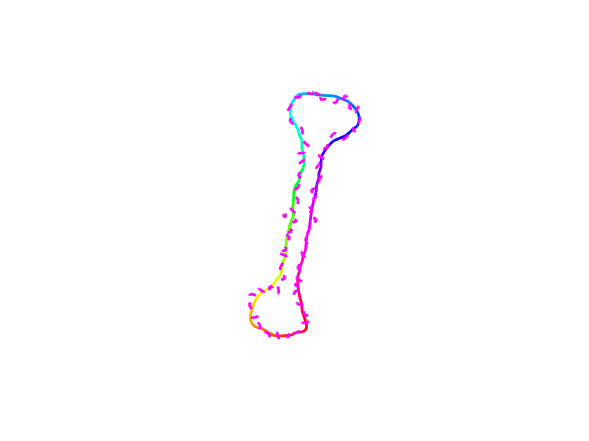}
\end{tabular}
\caption{Matching of curves with scale differences. The original noisy target curve is twice the size of the source. The first row shows the matching result obtained by length-normalization of the target (magenta curve). The last row is the registration obtained with the scale-invariant metric with simultaneous estimation of the rescaling parameter. The algorithm finds an optimal scaling of the target of 0.505 and the rescaled target is shown in magenta.} \label{fig:scale_invariance}
\end{figure}

Furthermore, even after choosing a ``good'' scale (or alternatively using a scale-invariant metric) the approach  of comparing curves at a fixed length might yield sub-optimal results.  In particular, in the presence of noise, rescaling the curves to constant length might artificially change their relative sizes.

Using scale-invariant metrics of the form \eqref{eq:sobolev_metric:scale_inv} overcomes both of these difficulties: the scale invariance of the metric makes the analysis independent of the choice of scale of the curves and it allows one to consider the induced Riemannian metric on the quotient space of curves modulo scalings by  optimizing over all rescalings of the target curve. Consequently this method automatically fits the optimal size of the target curve with respect to the relative size of the initial curve. 

This is shown in Figure \ref{fig:scale_invariance} where both strategies for dealing with scale variations are compared in a situation where the target curve's vertices are also corrupted by noise. As fidelity term, we use in this case a varifold metric with a linear function $\gamma$ which corresponds to the model of currents and was shown (c.f \cite{Charon2017}) to provide better robustness to such noise. As one can see in the first row of the figure curve length gives in that example a rather bad estimate for the rescaling factor and leads to a quite unnatural mapping where the usual cancellation effects of current fidelity metrics appear when trying to shrink the initial curve. In contrast, using a scale-invariant elastic metric (second row) allows one to jointly estimate a more sensible re-scaling parameter of the target shape (the variable $r$ in \eqref{eq:var_rigid_invariant}) together with a more natural path in the space of curves.

\subsection{Intrinsic vs extrinsic metric matching}
Another benefit of relaxing the constraint with varifold terms is that our new formulation can be more directly compared to another important class of shape space metrics and matching algorithms. We refer to those as ``extrinsic'' as they are usually related to the general model of shape spaces laid out by Grenander in \cite{Grenander1993}. In this model  distances and geodesics between two shapes are induced by a distance on a certain group of transformations of the entire embedding space that ``act'' on the shapes: in other words, the distance is quantified by the minimal amount of deformation necessary to map one shape to the other. 

\begin{figure}
\centering
\setlength{\tabcolsep}{.3em}
\begin{tabular}{ccccc}
\includegraphics[trim = 28mm 0mm 24mm 20mm ,clip,width=2.3cm]{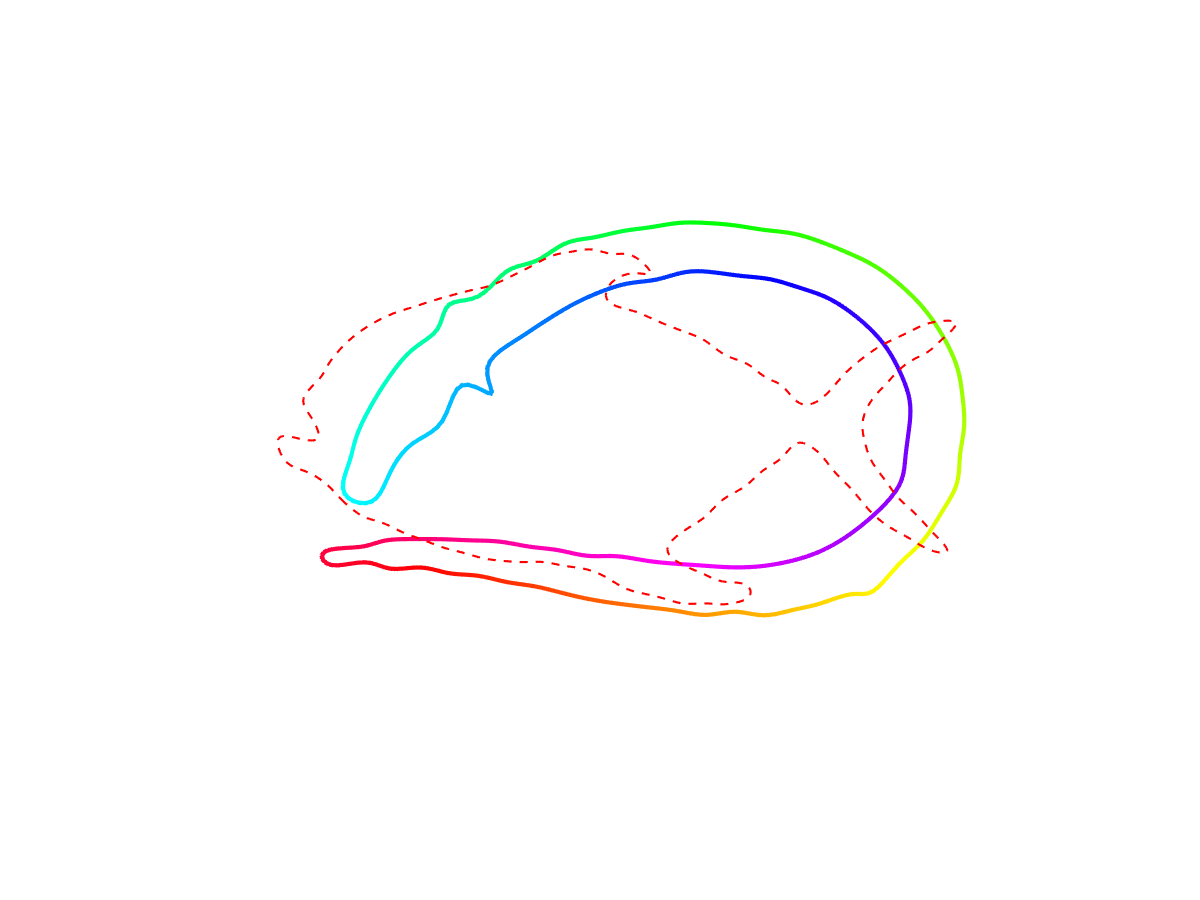} &
\includegraphics[trim = 28mm 0mm 24mm 20mm ,clip,width=2.3cm]{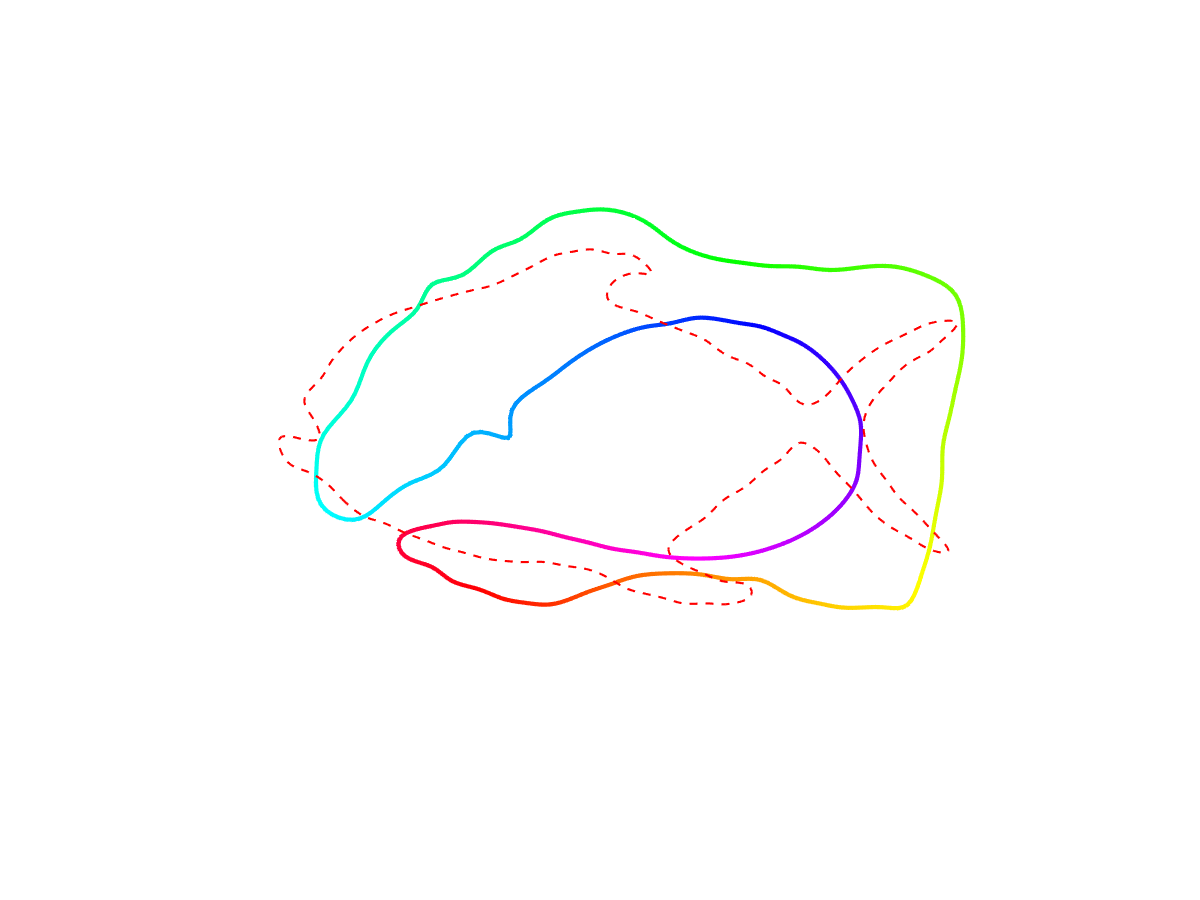} &
\includegraphics[trim = 28mm 0mm 24mm 20mm ,clip,width=2.3cm]{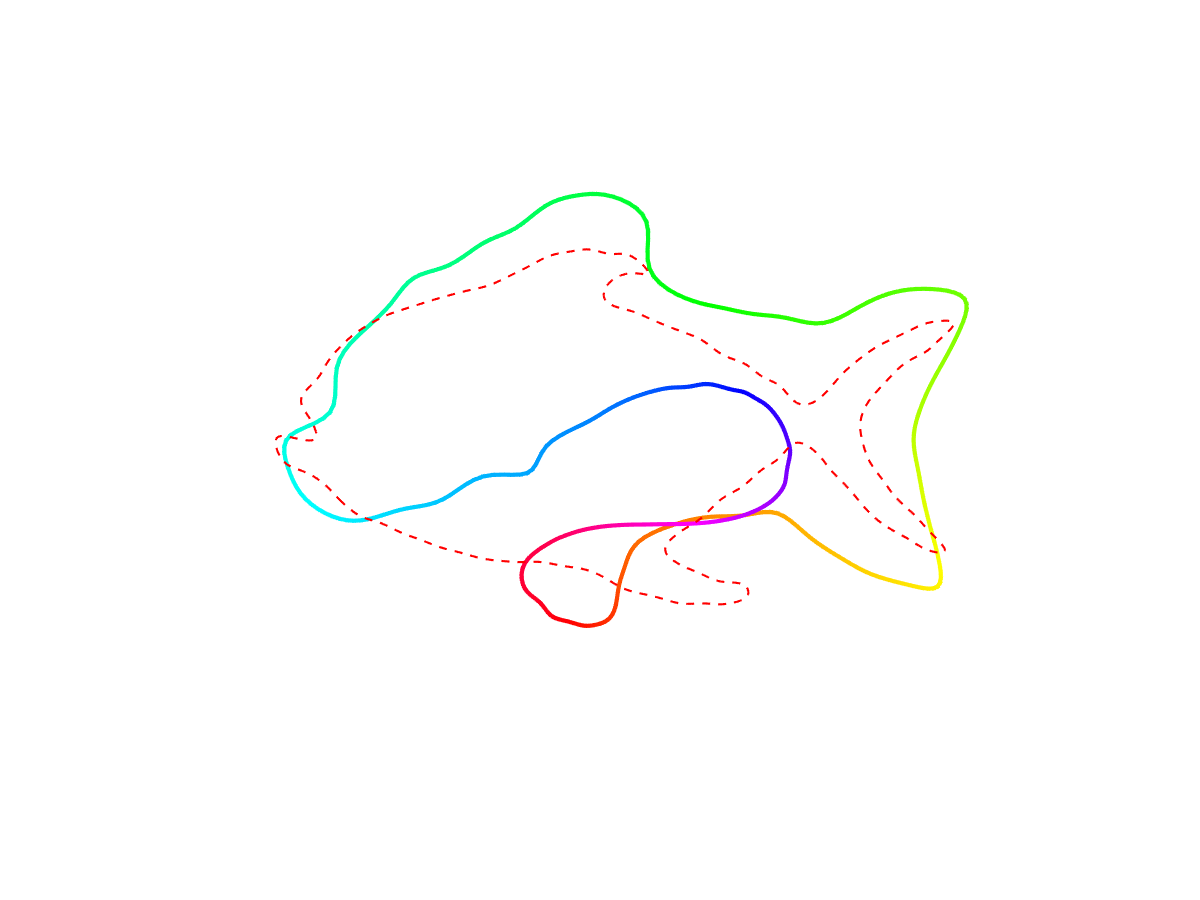} &
\includegraphics[trim = 28mm 0mm 24mm 20mm ,clip,width=2.3cm]{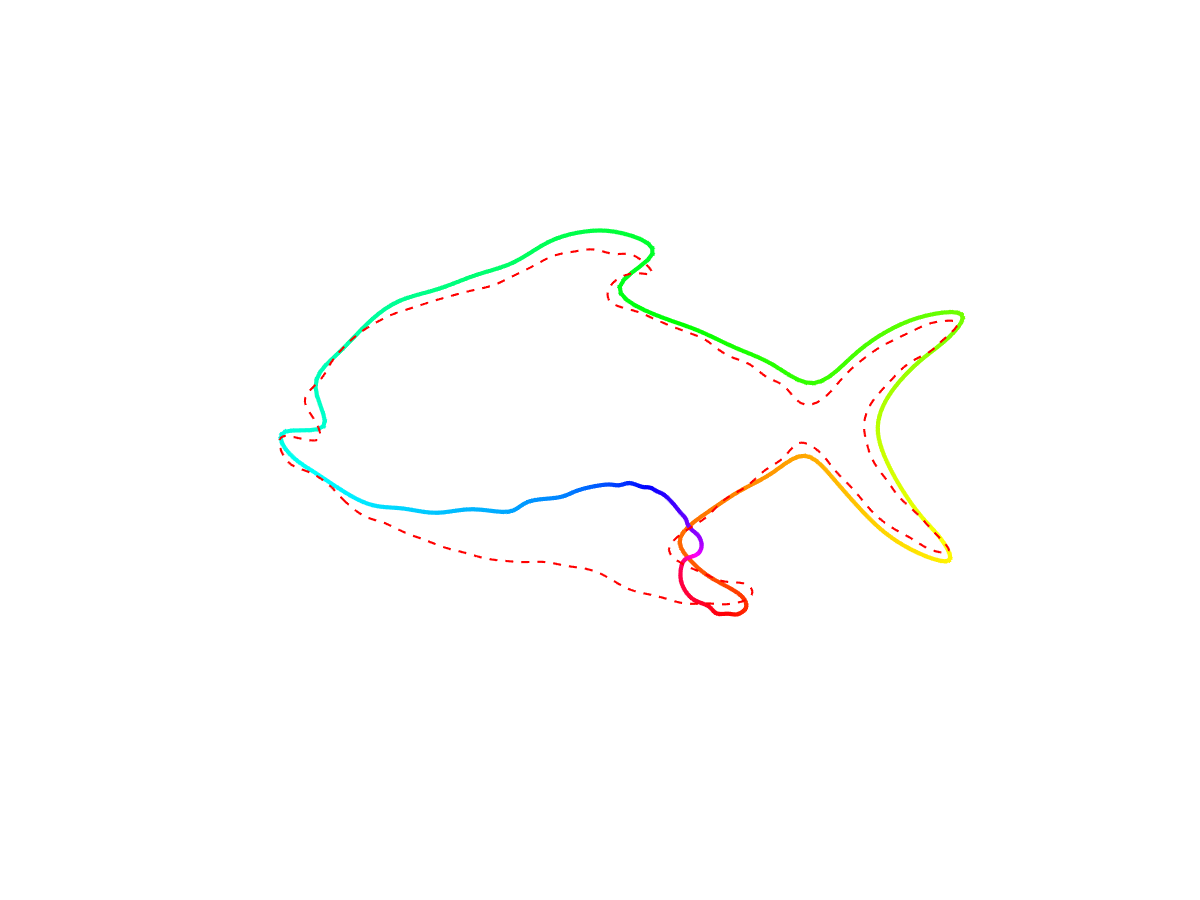}&
\includegraphics[trim = 28mm 0mm 24mm 20mm ,clip,width=2.3cm]{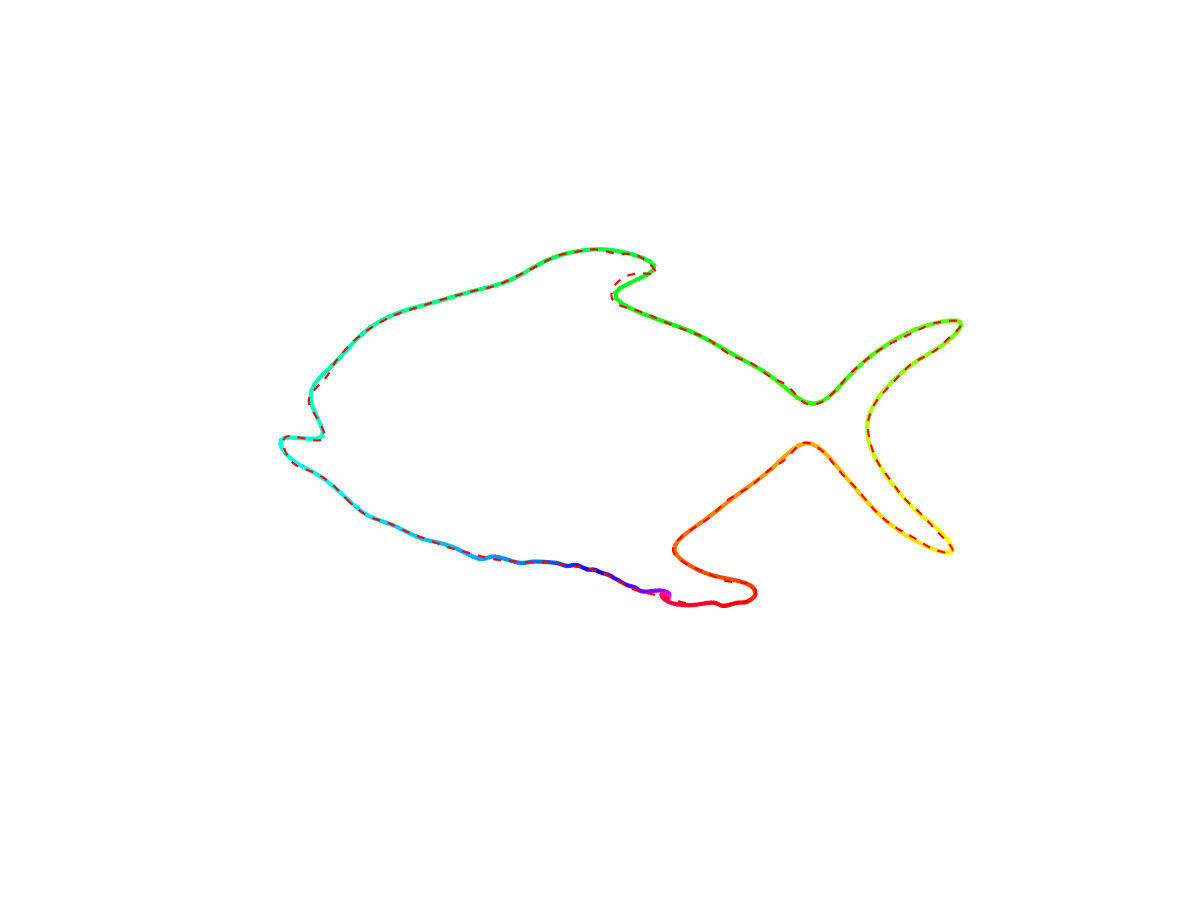}\\
t=0 & t=0.3 & t=0.6 & t=0.9 &t=1
\end{tabular}
\caption{An example of elastic Sobolev geodesic between two curves (the target is the red-dotted one). As opposed to extrinsic deformation models like LDDMM, self-intersections can be created in geodesic paths.} \label{fig:self_intersection}
\end{figure}

In most situations of interest, the groups in question are constructed as subgroups of $\on{Diff}(\R^d)$, the group of diffeomorphisms of the ambient space $\R^d$, equipped with a right-invariant metric. Multiple  models for such groups and metrics have been proposed. In this section we will focus on comparing our method with one of them: the Large Deformation Diffeomorphic Metric Mapping (LDDMM) framework originally introduced in \cite{Beg2005}. In the case of curves LDDMM inexact matching is typically formulated as the optimal control problem
\begin{equation}
    \label{eq:curve_matching_LDMMM}
    \inf_{v \in L^2([0,1],V)} \int_{0}^{1} \|v(t,\cdot)\|_V^2 dt + \lambda \|\mu_{c(1)} - \mu_{c_1}\|_{W^*}^2 
\end{equation}
on the time-dependent vector field $v \in L^2([0,1],V)$ where $V$ is a given reproducing kernel Hilbert space of smooth vector fields on $\R^d$, subject to the constraint $c(0)=c_0$ and $c_t = v(t,c(t))$. Note that the matching constraint is enforced again through a relaxation term based on the varifold metrics of Section \ref{sec:VarifoldDistance} (in fact \cite{Glaunes2008} uses metrics with $\gamma(u)=u$ while \cite{Charon2013} considers $\gamma(u)=u^2$ in applications). The essential difference between our formulation \eqref{eq:Energy_Vari} and \eqref{eq:curve_matching_LDMMM} is the fact that in \eqref{eq:curve_matching_LDMMM} the vector field is defined over the whole space $\R^d$ and its energy is measured by the global norm $\|\cdot\|_{V}$. 

This has a few important consequences. One key property of the LDDMM model is that it enforces the global transformation resulting from the flow of $v$ to be diffeomorphic. In particular, it will prevent any self-intersection from occurring along geodesics. In contrast, geodesics for the elastic Sobolev metrics of this paper lie in the space of immersions and, as illustrated in Figure \ref{fig:self_intersection}, self-intersections can appear in geodesics even if the initial and final curves are embedded curves.

\begin{figure}
\centering
\setlength{\tabcolsep}{.3em}
    \begin{tabular}{cccc}
\includegraphics[trim = 20mm 0mm 14mm 20mm ,clip,width=0.22\textwidth,height=0.25\textwidth]{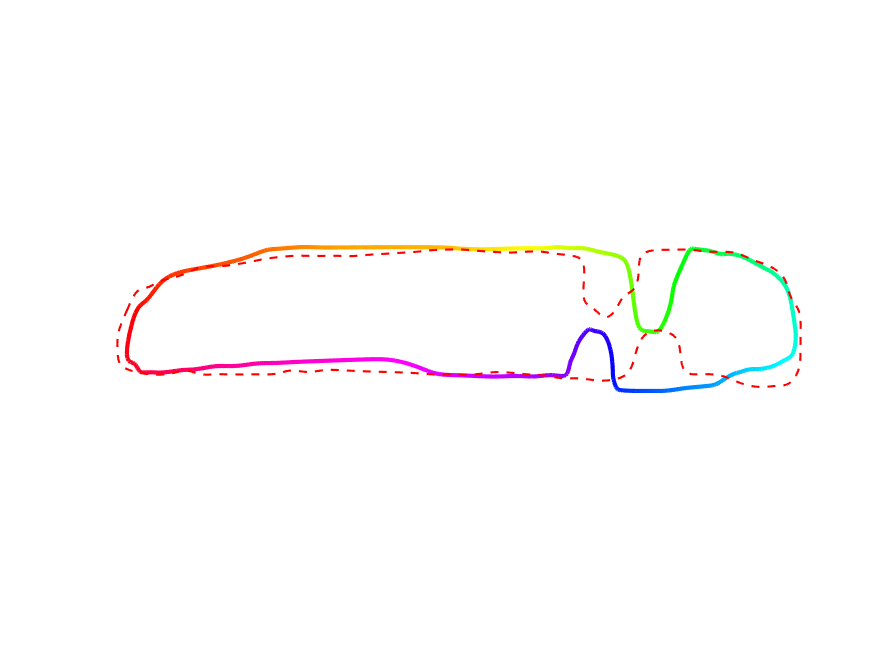} &
\includegraphics[trim = 20mm 0mm 14mm 20mm ,clip,width=0.22\textwidth,height=0.25\textwidth]{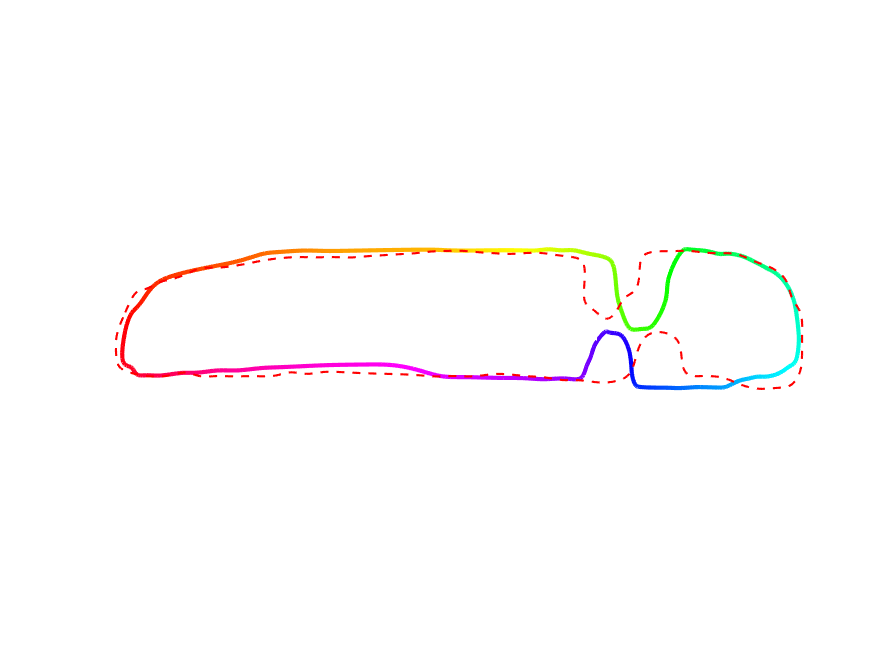} &
\includegraphics[trim = 20mm 0mm 14mm 20mm,clip,width=0.22\textwidth,height=0.25\textwidth]{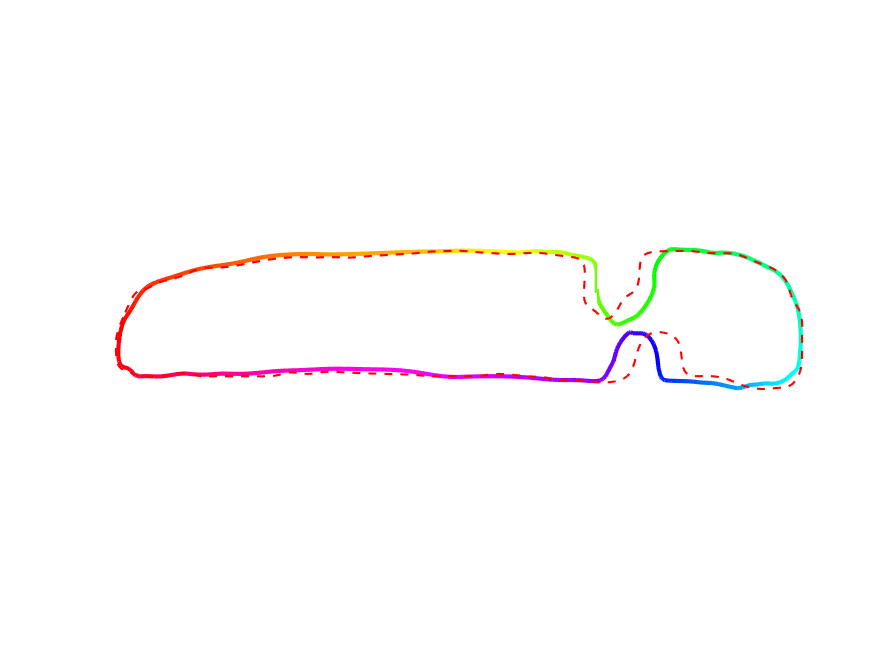} &
\includegraphics[trim = 20mm 0mm 14mm 20mm,clip,width=0.22\textwidth,height=0.25\textwidth]{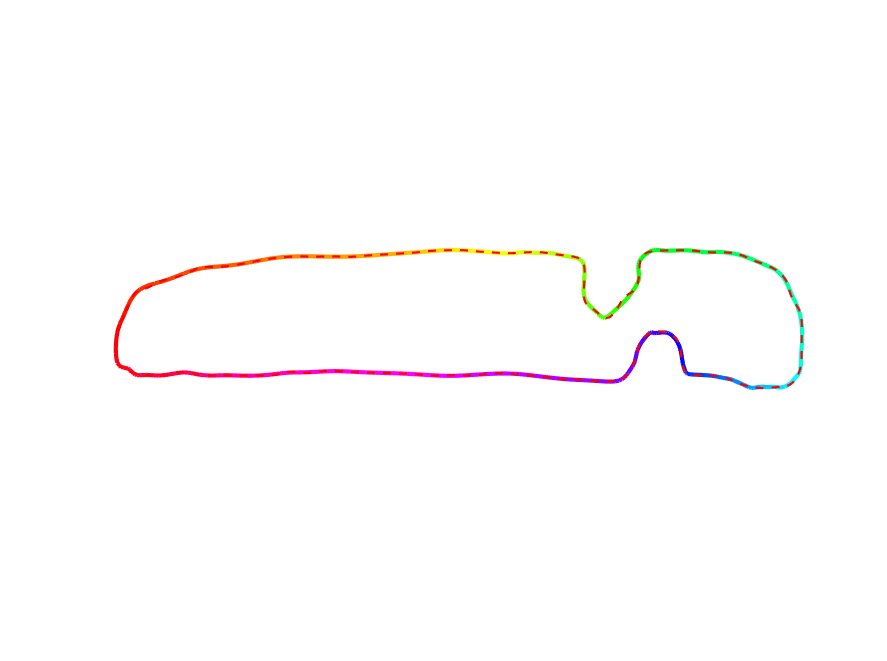}
\vspace{-1.5cm}\\
\includegraphics[trim = 20mm 0mm 14mm 20mm ,clip,width=0.22\textwidth,height=0.25\textwidth]{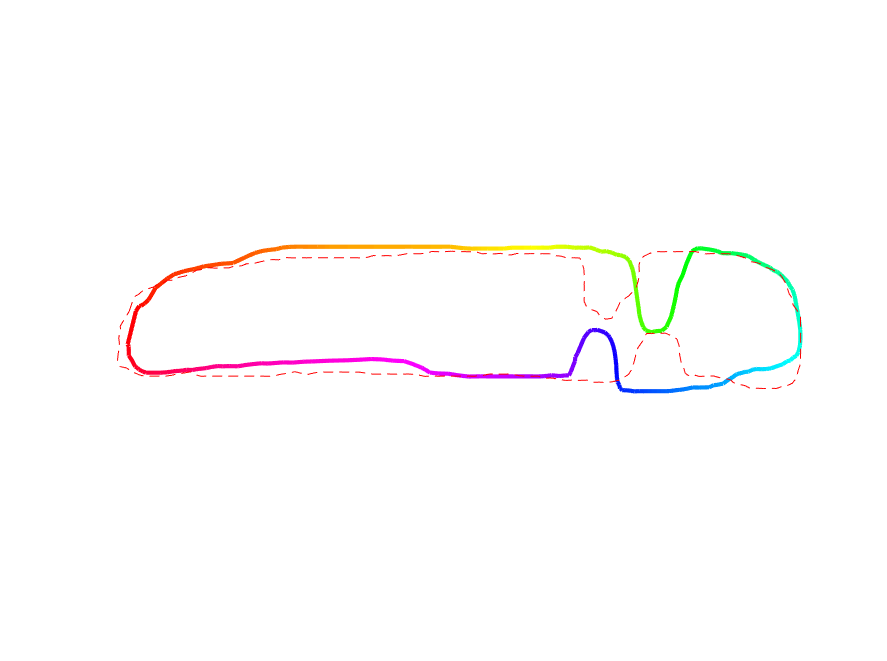} &
\includegraphics[trim = 20mm 0mm 14mm 20mm ,clip,width=0.22\textwidth,height=0.25\textwidth]{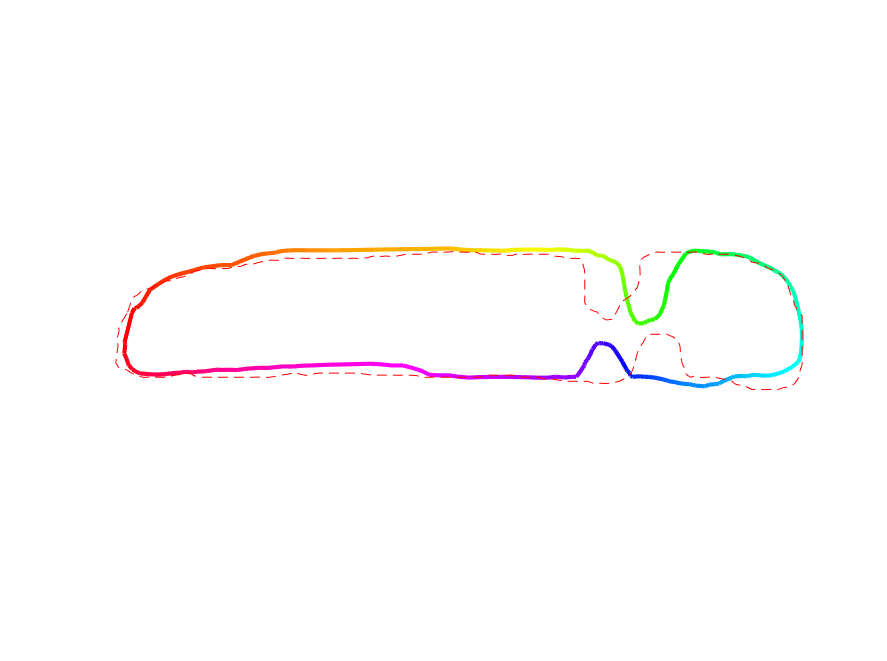} &
\includegraphics[trim = 20mm 0mm 14mm 20mm ,clip,width=0.22\textwidth,height=0.25\textwidth]{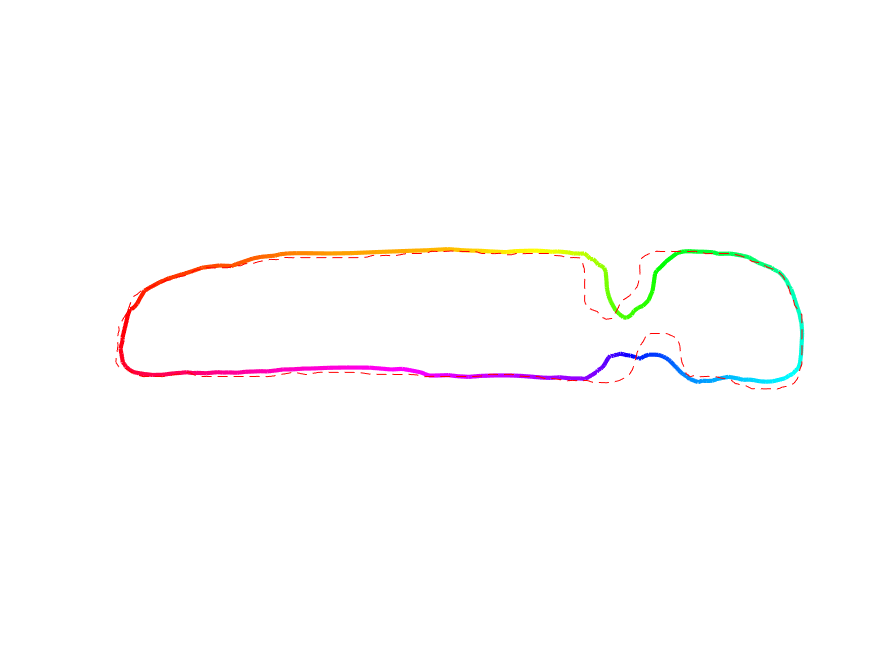} &
\includegraphics[trim = 20mm 0mm 14mm 20mm ,clip,width=0.22\textwidth,height=0.25\textwidth]{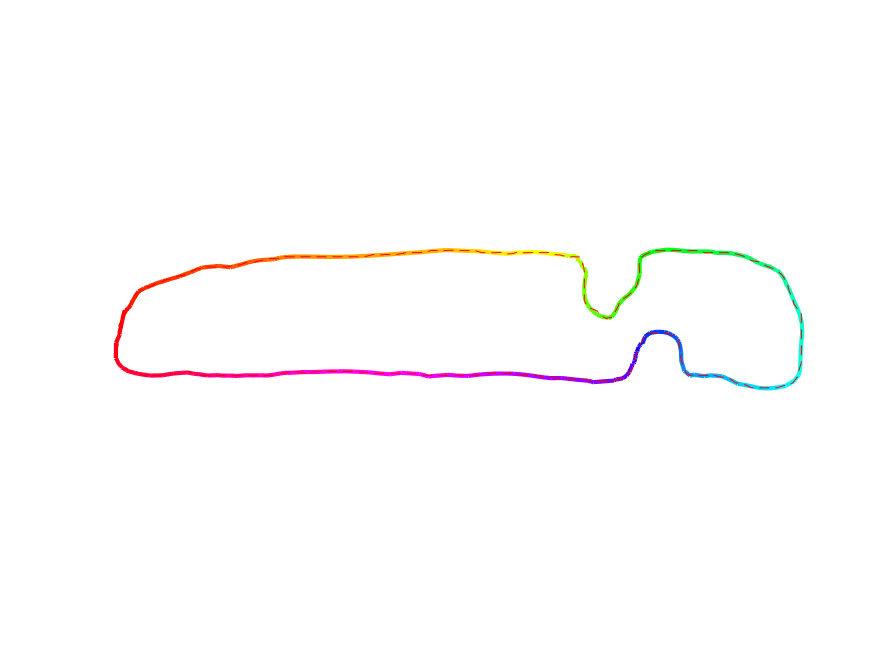}
\vspace{-1.5cm}
\\
$t=0$ & $t=0.3$ & $t=0.6$ & $t=1$ \vspace{-0cm}
    \end{tabular}
    \caption{Comparison of geodesics between curves with traveling bumps. First Line:Intrinsic $H^2$-metric; Second Line: LDDMM geodesic. In the LDDMM model one bump is successively flattened and recreated due to the high deformation cost of opposite displacements of close objects. In the intrinsic model the two bumps are merely transported.} \label{fig:traveling_bumps}
\end{figure}

\begin{figure}
\setlength{\tabcolsep}{.1em}
\begin{tabular}{cccc}
\includegraphics[trim = 25mm 10mm 20mm 10mm ,clip,width=3cm]{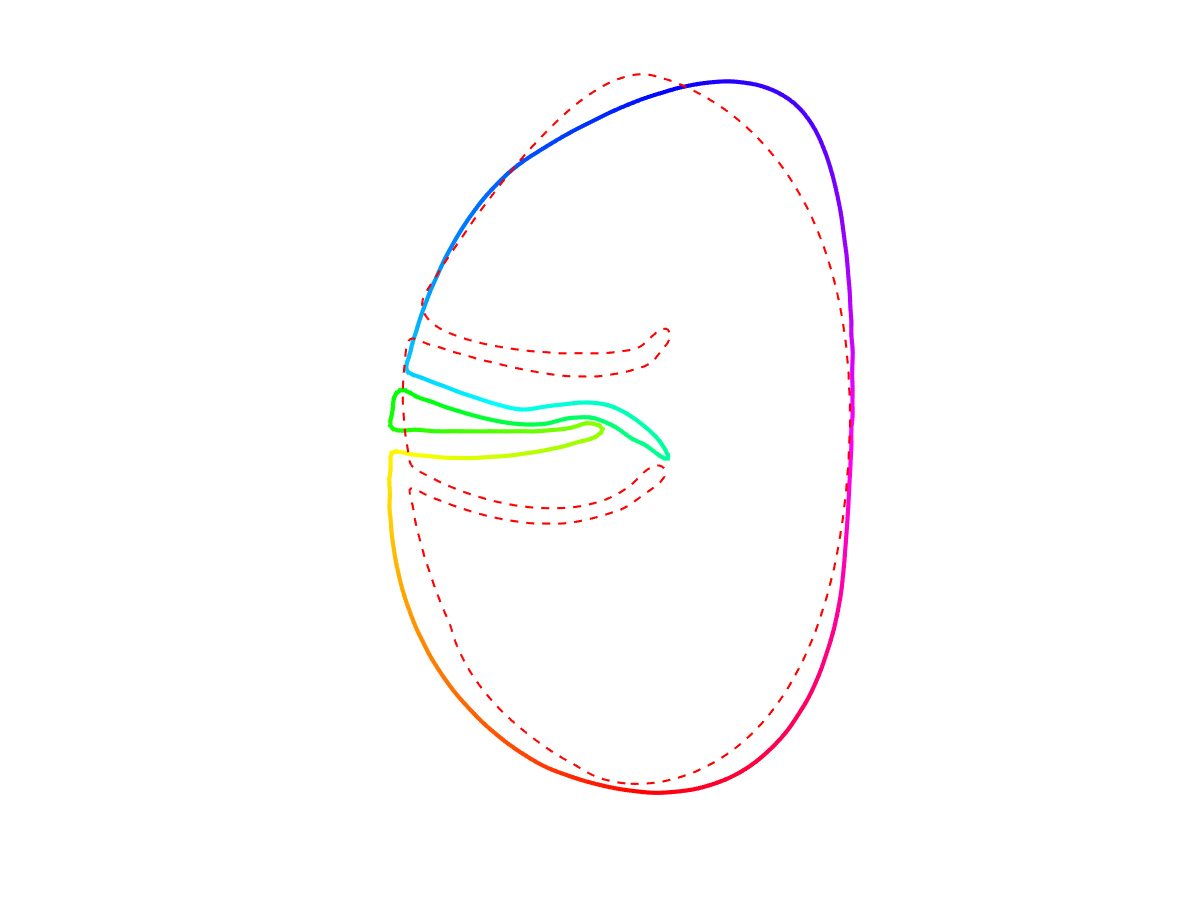} &
\includegraphics[trim = 25mm 10mm 20mm 10mm ,clip,width=3cm]{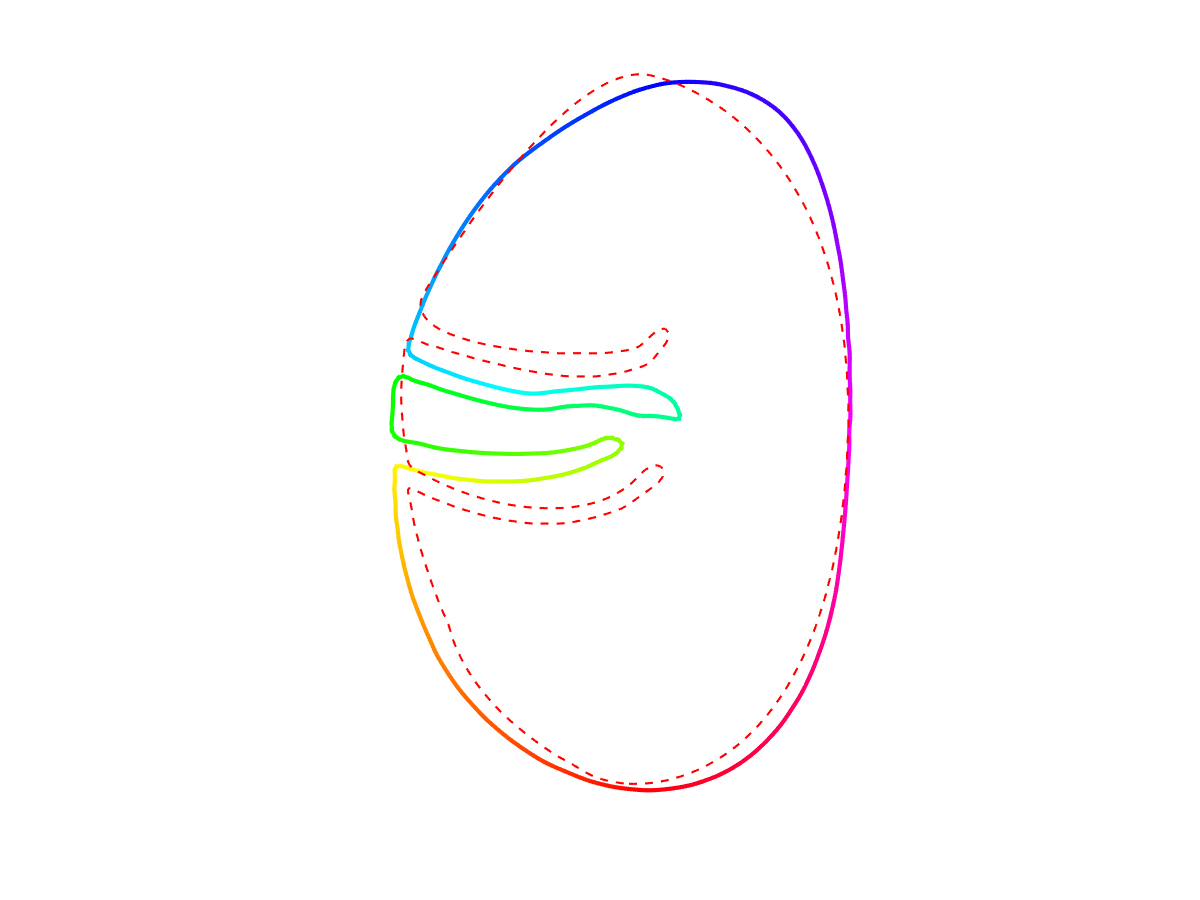} &
\includegraphics[trim = 25mm 10mm 20mm 10mm ,clip,width=3cm]{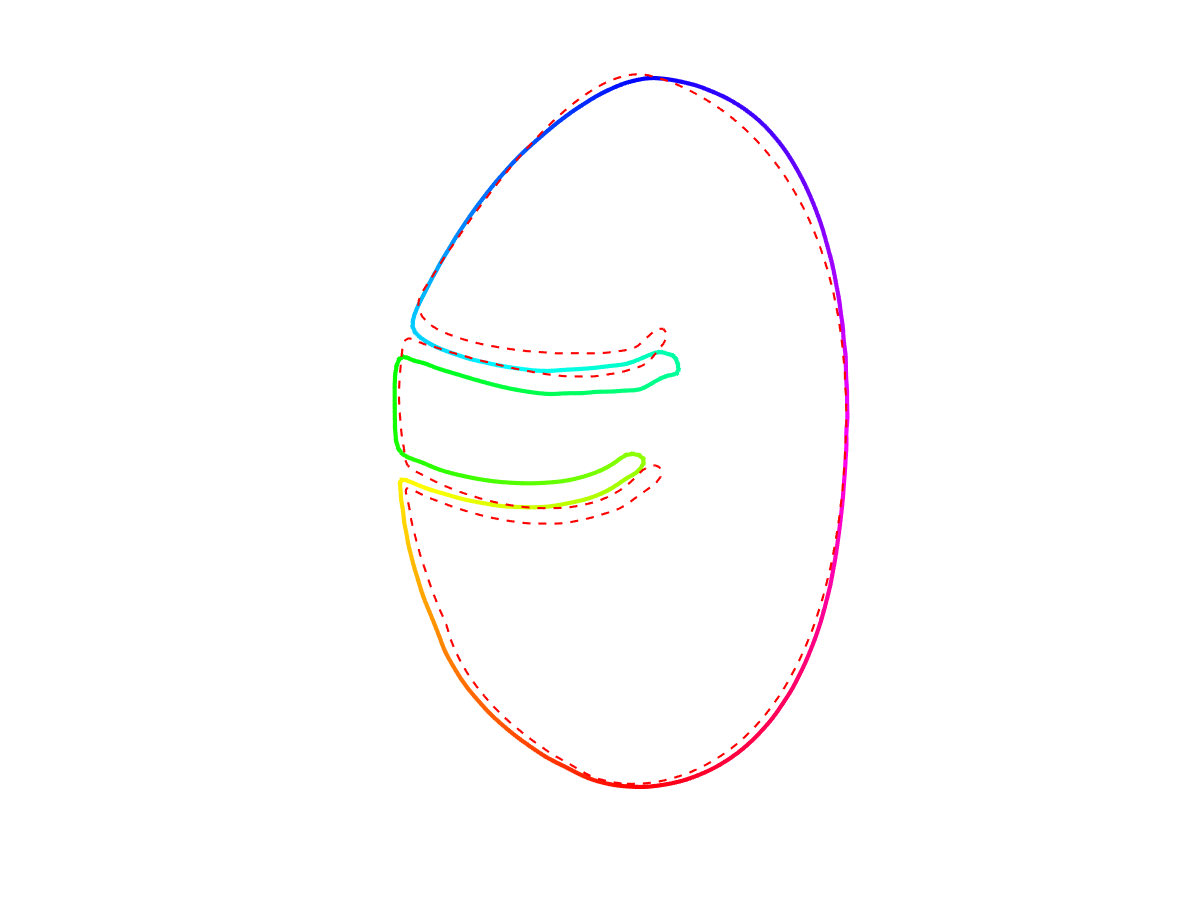} &
\includegraphics[trim = 25mm 10mm 20mm 10mm ,clip,width=3cm]{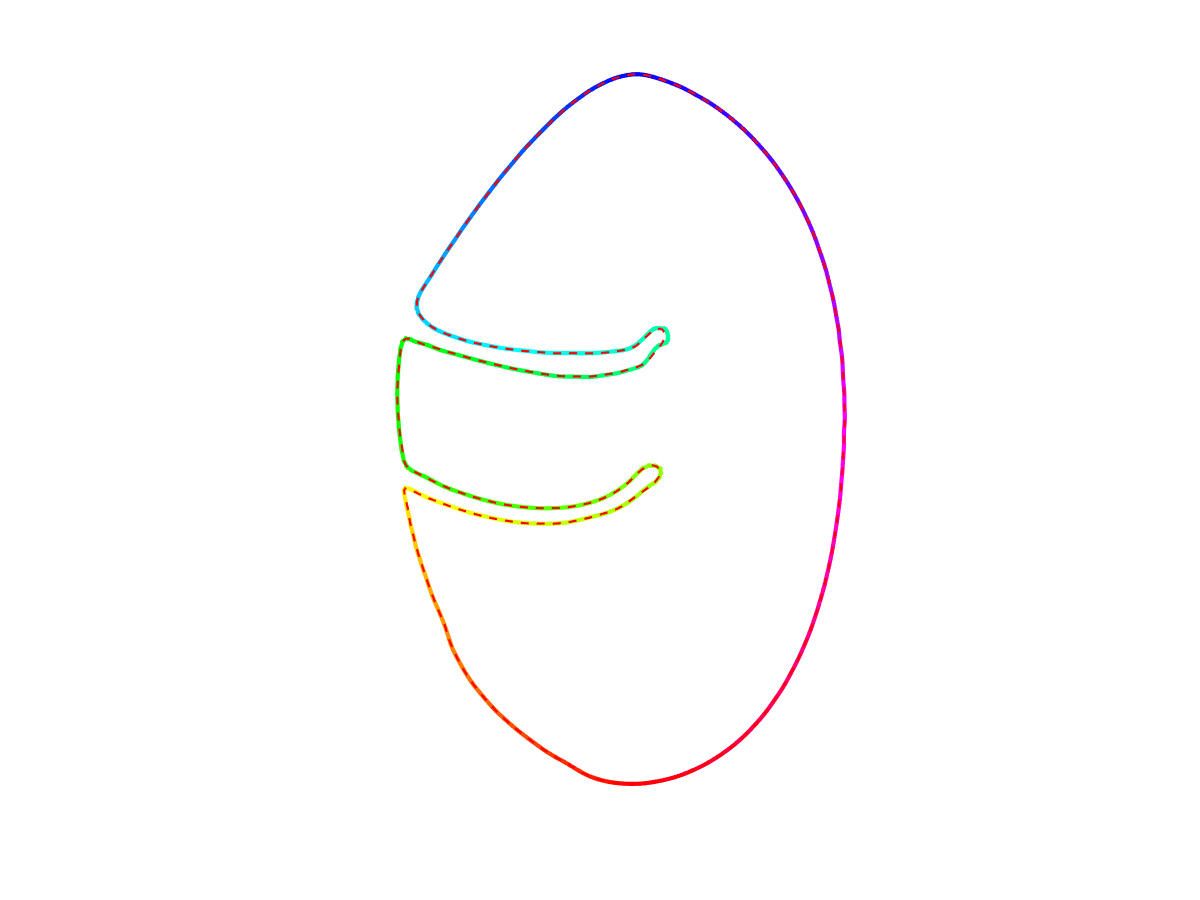}\\
$t=0$ & $t=0.3$ & $t=0.6$ & $t=1$
\end{tabular}
\caption{Estimated registration between the two curves with an intrinsic $H^2$ metric.} \label{fig:sulci_h2}
\end{figure}

\begin{figure}
\setlength{\tabcolsep}{.1em}
\begin{tabular}{cccc}
\includegraphics[trim = 25mm 10mm 20mm 10mm ,clip,width=3cm]{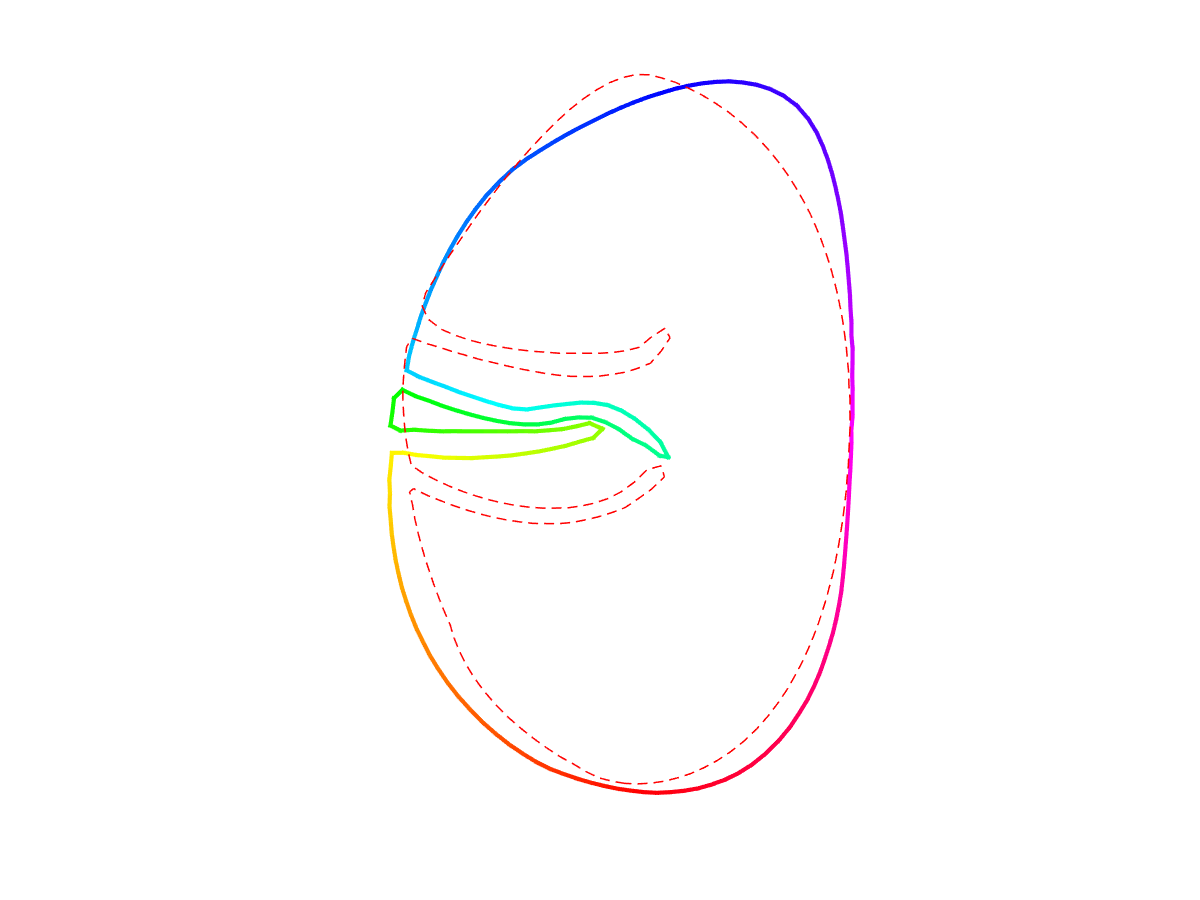} &
\includegraphics[trim = 25mm 10mm 20mm 10mm ,clip,width=3cm]{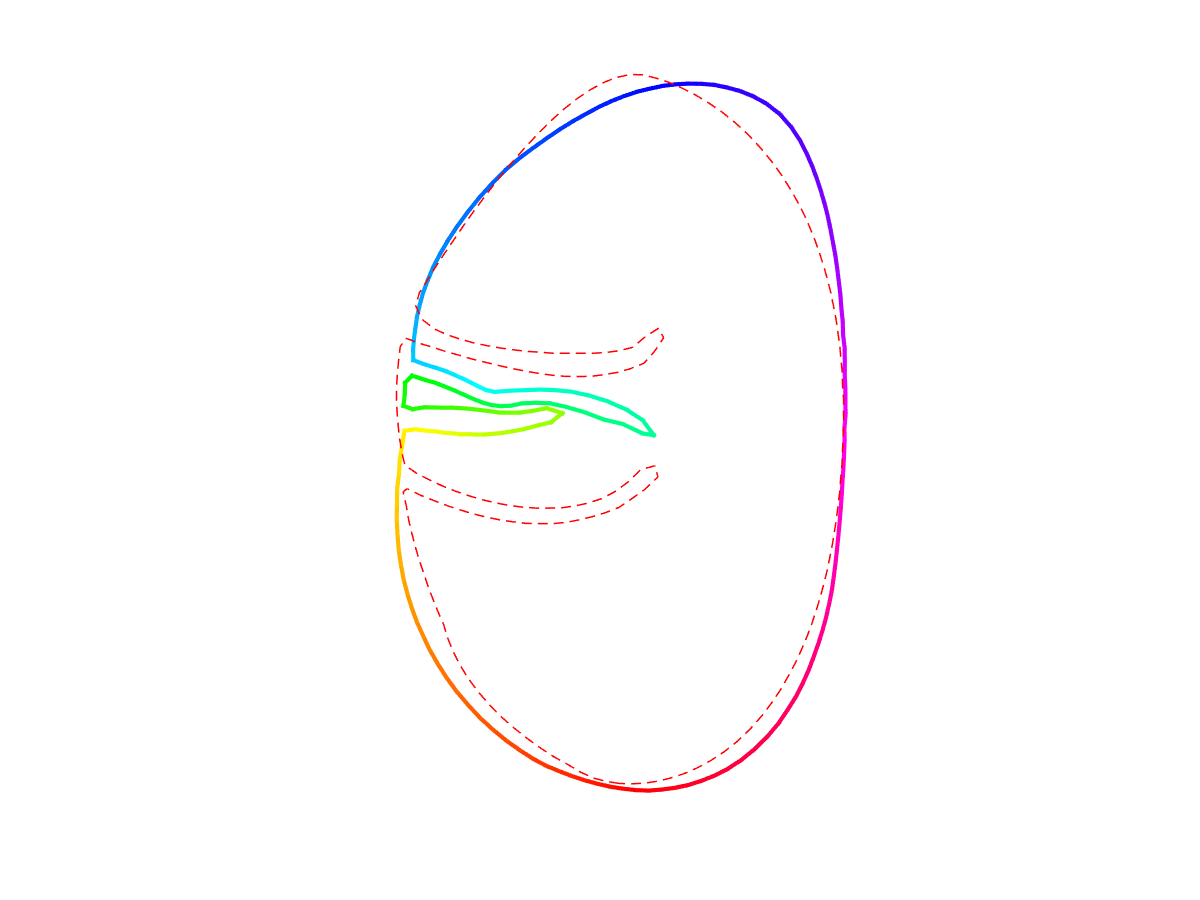} &
\includegraphics[trim = 25mm 10mm 20mm 10mm ,clip,width=3cm]{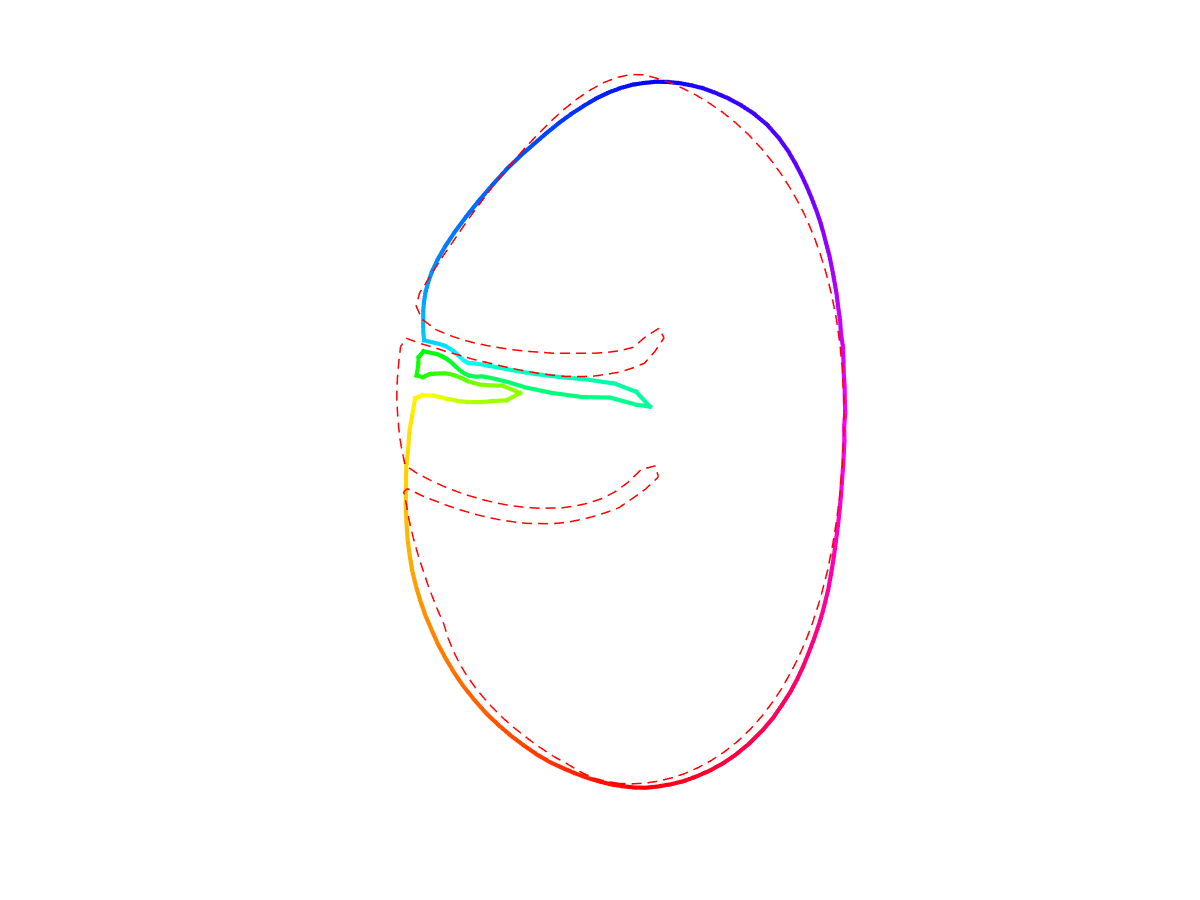} &
\includegraphics[trim = 25mm 10mm 20mm 10mm ,clip,width=3cm]{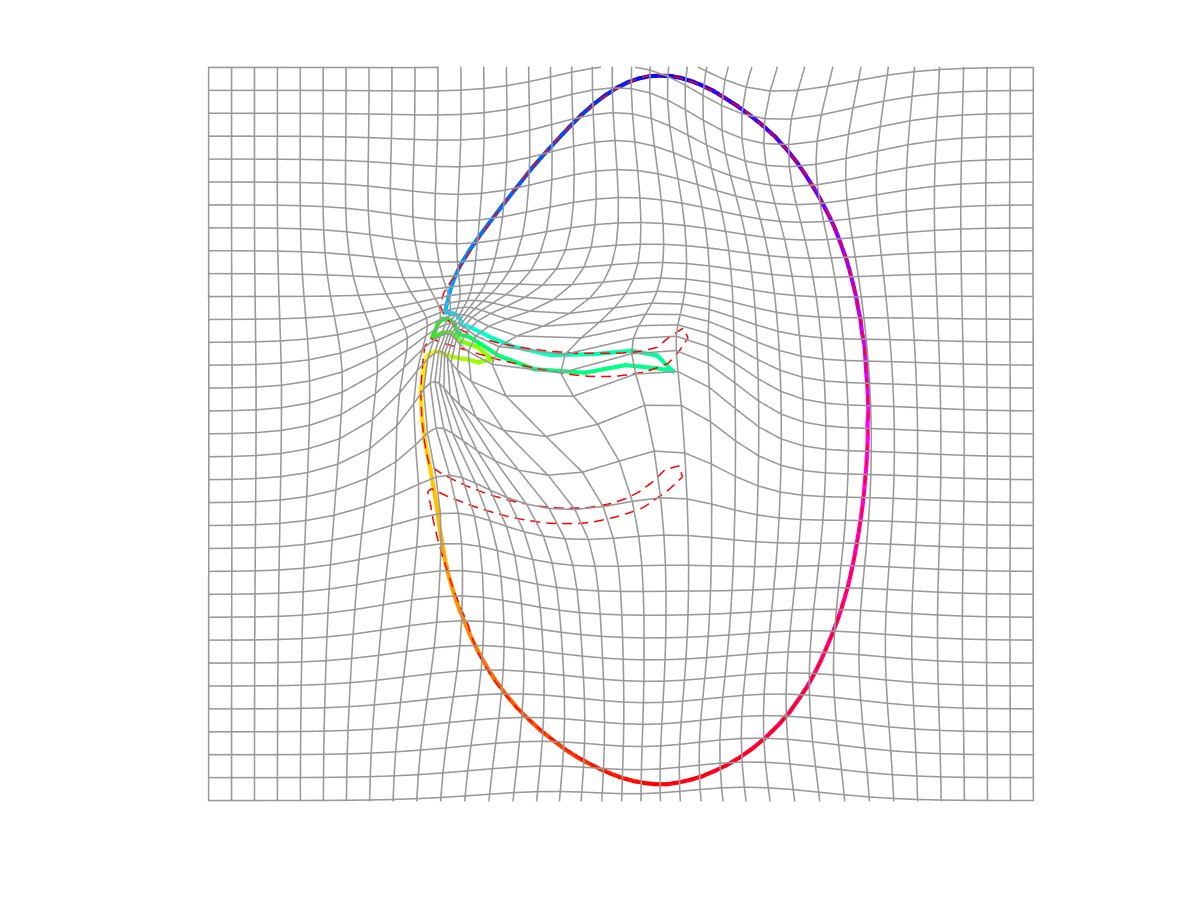}\\
$t=0$ & $t=0.3$ & $t=0.6$ & $t=1$ \\
\includegraphics[trim = 25mm 10mm 20mm 10mm ,clip,width=3cm]{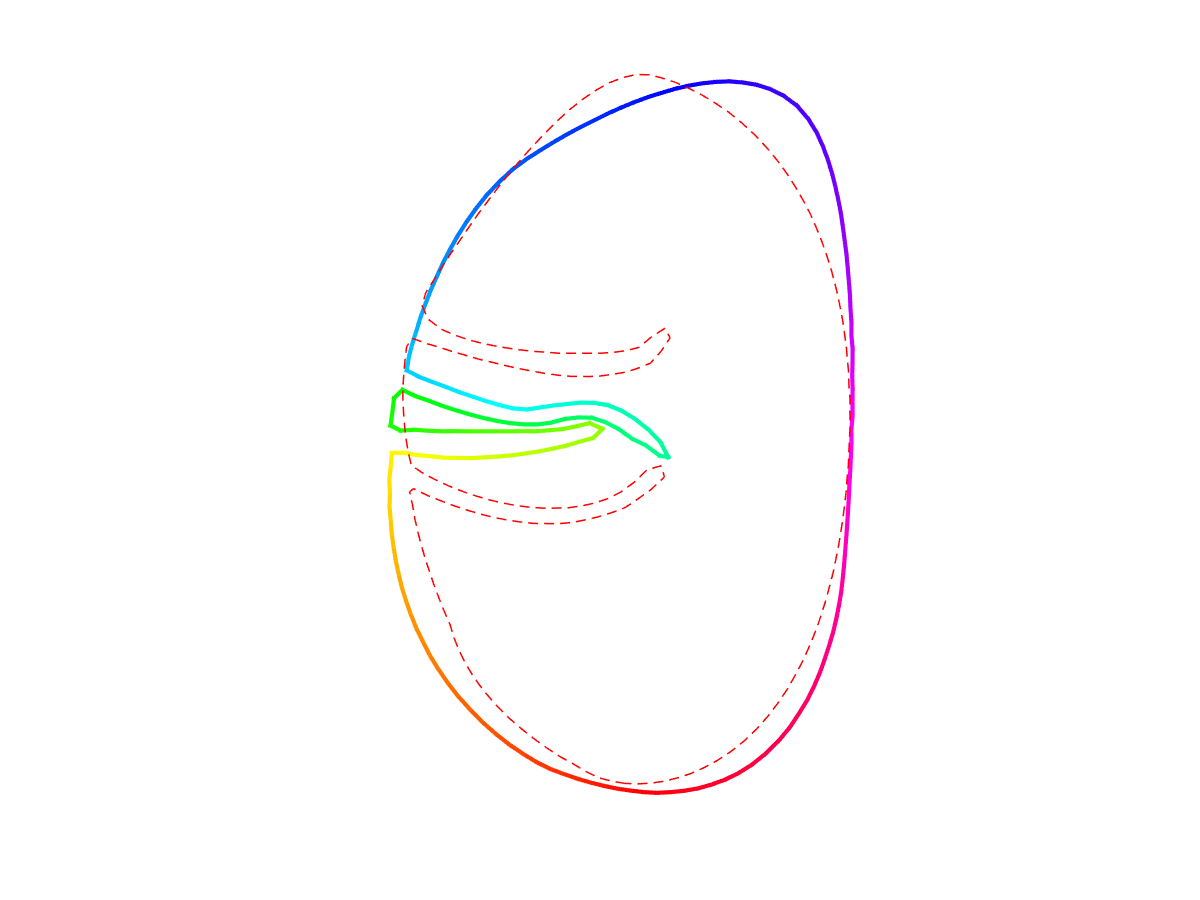} &
\includegraphics[trim = 25mm 10mm 20mm 10mm ,clip,width=3cm]{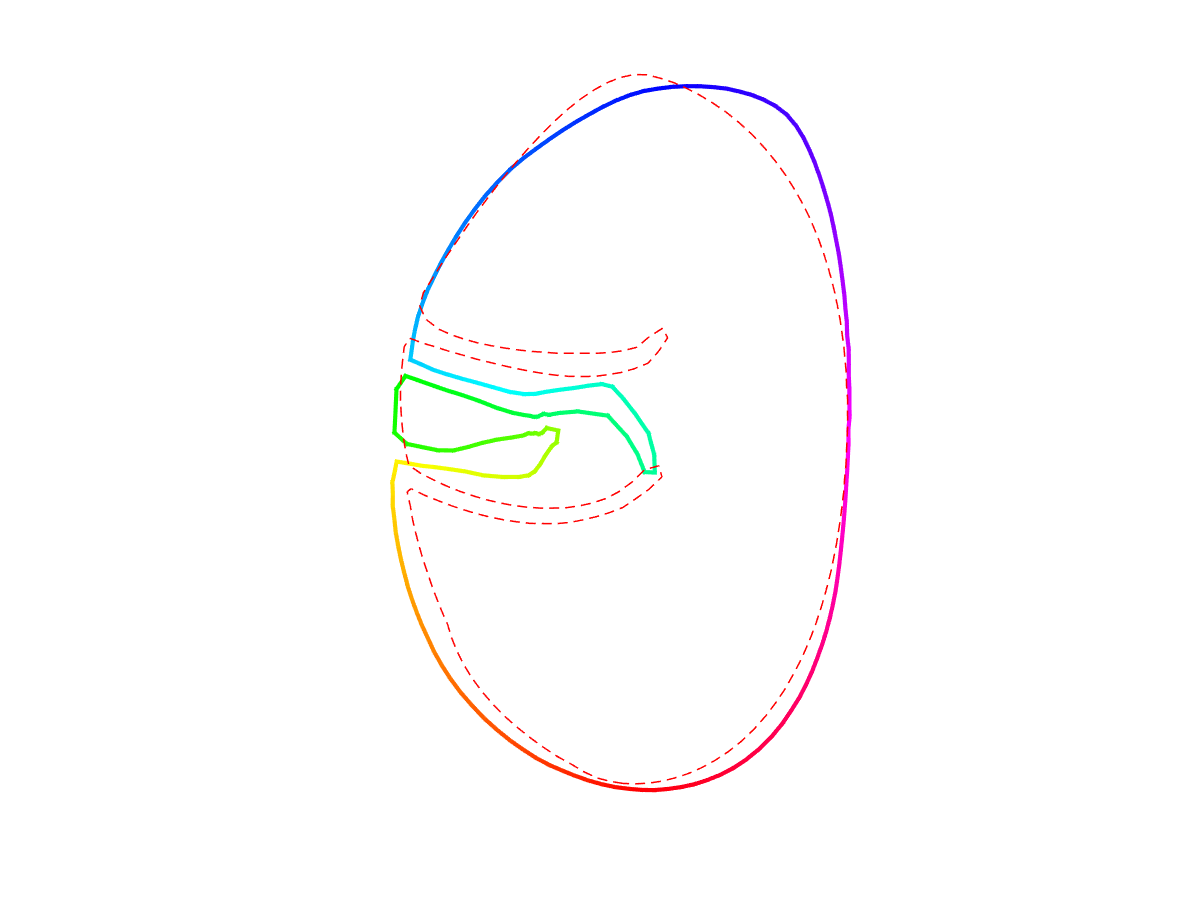} &
\includegraphics[trim = 25mm 10mm 20mm 10mm ,clip,width=3cm]{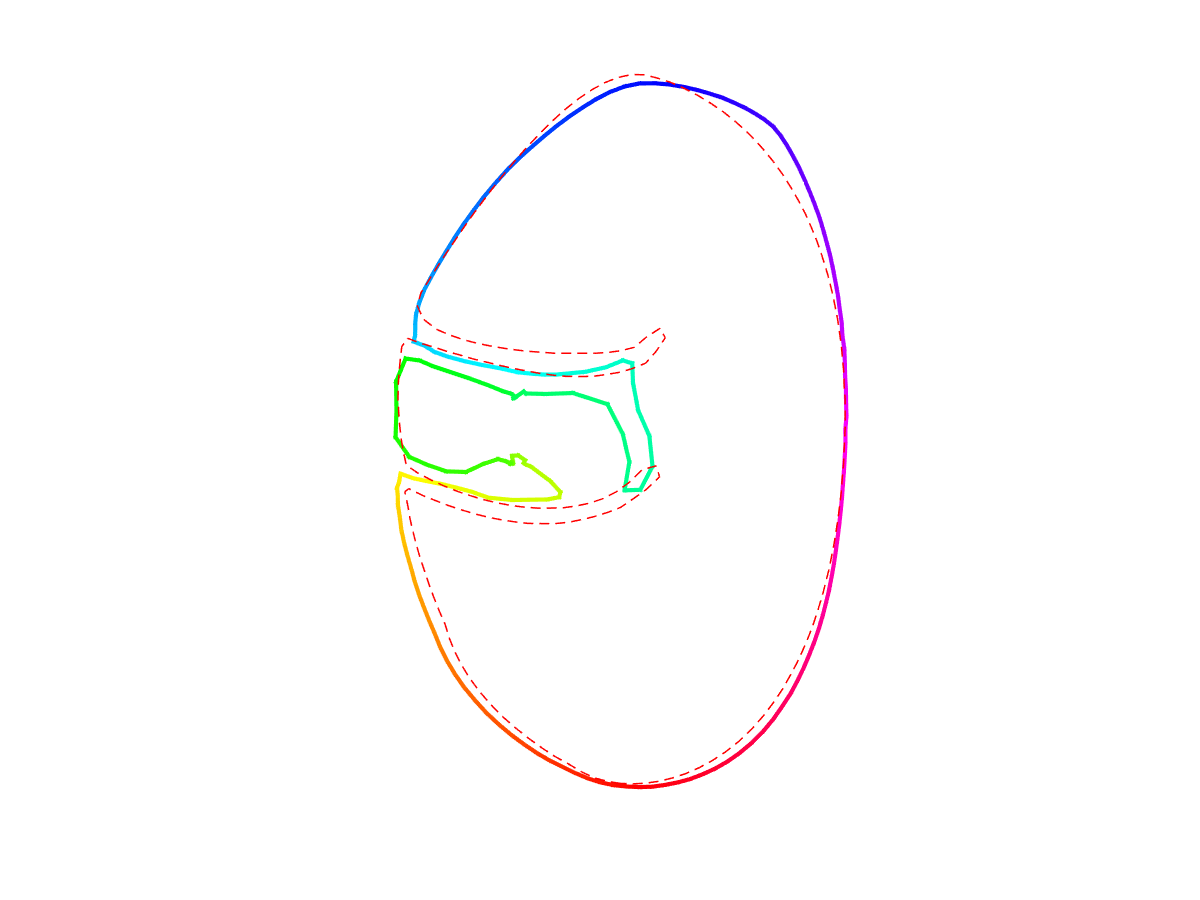} &
\includegraphics[trim = 25mm 10mm 20mm 10mm ,clip,width=3cm]{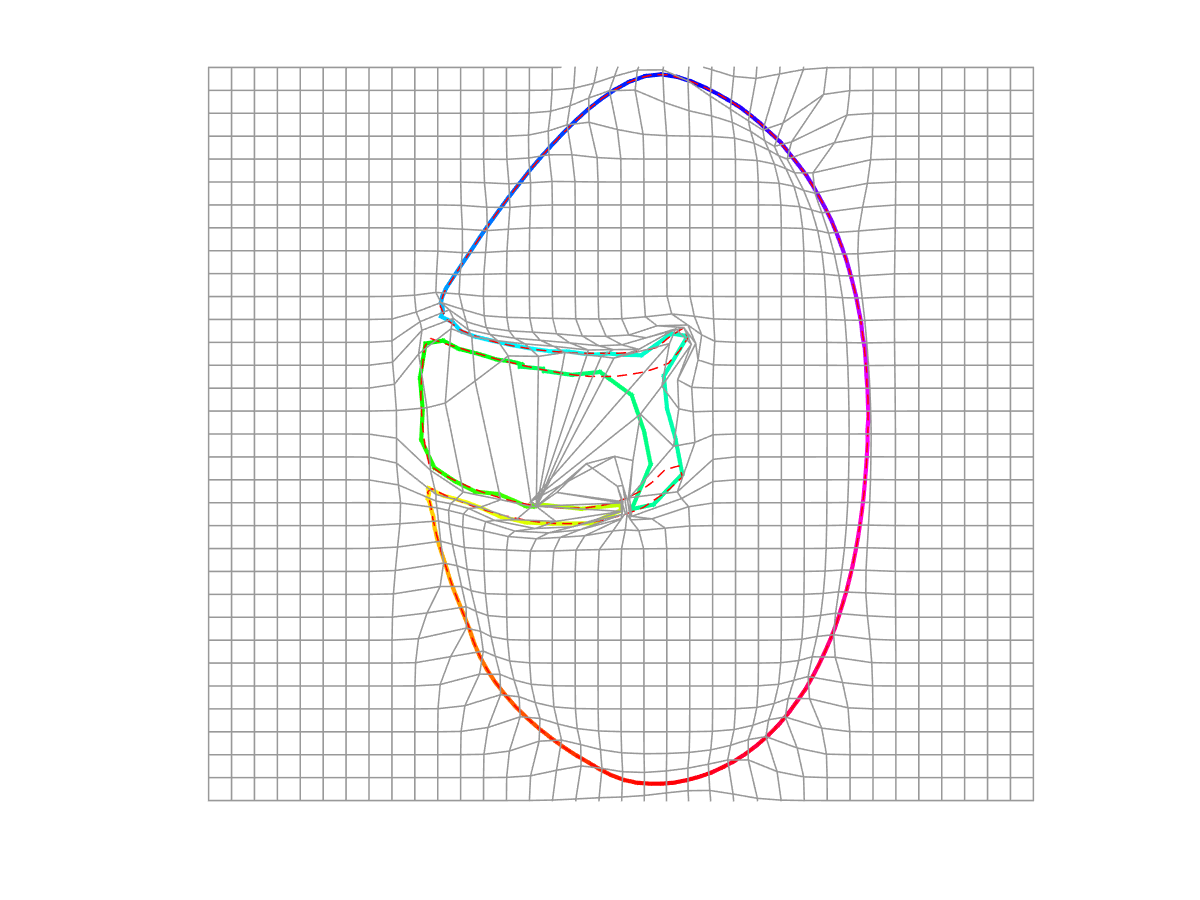}
\end{tabular}
\caption{Registration using LDDMM. On the first row, using a Gaussian deformation kernel of width $0.2$. On the second row, with a Gaussian kernel of width $0.05$. Note that if reducing the kernel size allows theoretically to recover finer scale deformations, it is well-known that the LDDMM registration problem then becomes highly sensitive to bad local minima such as the solution shown in the second row.} \label{fig:sulci_lddmm}
\end{figure}

On the other hand estimating a global diffeomorphic transformations as given by the LDDMM model may prove a particularly difficult or undesirable constraint in certain situations. It is most notably the case when thin or closely located structures have to be displaced or stretched apart. We illustrate such a phenomenon in the example of Figure \ref{fig:traveling_bumps}. The motion of bumps on two opposite sides of a curve is estimated in fundamentally different ways by the intrinsic model for which the two bumps are simply displaced along the curve and by LDDMM where one bump is successively flattened and recreated due to the high deformation cost of opposite displacements of close by objects. Along the same lines, Figures \ref{fig:sulci_h2} and \ref{fig:sulci_lddmm} show another comparison in which two ``sulci'' have to be moved apart. Since this is again a costly deformation in the LDDMM framework, it is easily prone to reach unnatural solutions if the deformation kernel is too large or to lead to even more unnatural local minima of the functional \eqref{eq:curve_matching_LDMMM} for small deformation kernels. We also point out that similar issues are discussed quite extensively in the recent work \cite{Younes2017}, which in addition introduces a hybrid model combining a global LDDMM deformation cost with intrinsic $H^1$ penalties.      

\subsection{Time comparison}
As a last set of experiments, we take a closer look at the computational complexity for both of the previous registration models. As detailed earlier in Section \ref{sec:numerics}, the bulk of the computations for the proposed elastic metric approach at each optimization step is 1) the evaluation of splines and their derivatives to compute the elastic energy and its gradient for which the required number of basic operations is linear in the total number $N_\theta N_t$ of spline control points and 2) the evaluation of the varifold fidelity term which requires $O(N_\theta^2)$ evaluations of the kernels $\rho$ and $\gamma$. Standard LDDMM algorithms, in comparison, usually solve \eqref{eq:curve_matching_LDMMM} with a forward-backward shooting procedure involving, at each iteration of the optimization, the integration of Hamiltonian systems of equations with interacting particles, which will typically require of the order of $O(N^2 T)$ evaluations of the kernel of the RKHS $V$, with $N$ being the number of vertices of the curve and $T$ the number of time steps for the chosen numerical integrator.

\begin{table}
\centering
\setlength\extrarowheight{3pt}
\begin{tabular}{lcccc}
\toprule
  & Total run time & Avg. $n\textsuperscript{o}$ of iter./registration & Avg. time/iter. \\
\midrule
Elastic & 66 mins & 479 & 0.02s  \\
LDDMM & 93 mins & 91 & 0.16s \\ 
\bottomrule
\end{tabular}
\vspace{5pt}
\caption{Computational time comparison between the proposed algorithm and LDDMM for 380 pairwise registrations of curves with 100 vertices each. We report the total run time, the average number of iterations for the optimization methods and the average time of one iteration.}
\label{tab:time}
\end{table}

We illustrate this comparison empirically by running the two registration algorithms on a subset of 20 closed curves from the Surrey fish dataset and estimating all the 380 pairwise matchings for each method. Note that, although the two optimization problems and part of their underlying parameters differ, the quality of the registration is in both cases measured by the varifold metric. In this experiment, we choose the same scale parameters for the fidelity term and adapt the other parameters relative to the deformation metrics to lead comparable convergence properties and registration accuracy for both methods. All curves are set to $N=100$ vertices which are also the control points used in the spline representation (thus $N_\theta = N$), and we take $T=10$ time steps for LDDMM as well as $N_t=10$ time control points for splines in our proposed algorithm. The implementation of curve matching LDDMM is the one of \cite{Charon2017,fshapesTk} with the optimization routine given by the same limited memory BFGS algorithm from the HANSO library that is used in this work.

The results are reported in Table \ref{tab:time}. There are a few remarks to be made. First, on average, the time for a single iteration of the optimization procedure is significantly lower with the approach of this paper compared to LDDMM, which is consistent with the previous discussion on the theoretical complexity for the computation of deformation energies and gradients in both cases. Second, still on average, the LDDMM algorithm requires a priori less iterations for BFGS to reach convergence, with the same stopping criterion being used. This is likely due, on the one hand, to the fact that the optimization in LDDMM is performed over the deformation's initial momenta as opposed to the full path of spline parameters in our approach, thus reducing the size of the problem. On the other hand, it is also important to point out that this may be in part due to BFGS occasionally converging to irrelevant local minima (in very few number of iterations) in the case of LDDMM. In this precise experiment, this happens for about 30 registration cases in which the residual varifold cost remains very high at the end of the minimization. In contrast, the convergence seems much more consistent in the case of the elastic method as the total number of iterations and final energies do not vary as significantly from one registration to another.

As additional future comparison, it will be interesting to investigate the influence of the choice of $N_\theta$ or $N$ on the computational time and convergence properties for the two models. We postulate that the use of splines represented by their control points instead of directly vertices could allow $N_\theta$ to be in practice much smaller than $N$ while still providing consistent registration results for smooth curves.   
\section*{Acknowledgments}
We would like to thank Philipp Harms, Eric Klassen, Sebastian Kurtek, Peter Michor, Tom Needham, Anuj Srivastava and the Shape Group at FSU for helpful comments and discussions. Nicolas Charon is supported by the National Science Foundation under Grant No 1819131.

\appendix
\section{Derivatives of the energy functional}
\label{sec:Derivatives}
In this appendix we list the derivatives of the energy 
functional~\eqref{eq:EnergyFunctional} and varifold distance \eqref{eq:metric_W_curves}. The first derivative of the energy
\label{sec:AppendixEnergyDerivatives}
\begin{align*}
dE_c(k) = \int_0^1 \int_0^{2\pi} &\!\!t_1 \ip{c'}{k'} + t_2 \left( \ip{c''}{k'} + \ip{c'}{k''} \right)  
+t_3 \ip{\dot{c}'^{\bot}}{k'}
+ t_4 \ip{\dot{c}}{\dot{k}} + t_5 \ip{\dot{c}'}{\dot{k}'} \\
&+t_6\ip{\dot c'^{\top}}{\dot k'}+t_7\ip{\dot c'^{\bot}}{\dot k'}
+ t_8 ( \ip{\dot{c}''}{\dot{k}'} + \ip{\dot{c}'}{\dot{k}''} ) + t_9 \ip{\dot{c}''}{\dot{k}''}  \\
& +  d\ell_c(k) \bigg(  a_0'(\ell) |c'|\ip{\dot{c}}{\dot{c}} + a_1'(\ell)\frac{\ip{\dot{c}'^\top}{\dot{c}'^\top}}{|c'|} + b_1'(\ell) \frac{\ip{\dot{c}'^\perp}{\dot{c}'^\perp}}{|c'|} \\ &+ \left. a_2'(\ell) \left( \frac{\ip{c'}{c''}^2}{|c'|} + \frac{\ip{\dot{c}''}{\dot{c}''}}{|c'|^3} - \frac{2\ip{c'}{c'}\ip{\dot{c}''}{\dot{c}'}}{|c'|^5}  \right) \right) \ud \theta \ud t\,,
\end{align*}
with
\begin{align*}
t_1 &=  \frac{a_0}{|c'|} \ip{\dot{c}}{\dot{c}} -
\frac{a_1}{|c'|^3} \ip{\dot{c}'^\top}{\dot{c}'^\top}
-\frac{b_1}{|c'|^3} \ip{\dot{c}'^\bot}{\dot{c}'^\bot}
 - 7\frac{a_2}{|c'|^9} \ip{c'}{c''}^2 \ip{\dot{c}'}{\dot{c}'} \\
 &\qquad
 + 10\frac{a_2}{|c'|^7} \ip{c'}{c''} \ip{\dot{c}'}{\dot{c}''}
 - 3\frac{a_2}{|c'|^5} \ip{\dot{c}''}{\dot{c}''} \,,\\
t_2 &=  2\frac{a_2}{|c'|^7} \ip{c'}{c''} \ip{\dot{c}'}{\dot{c}'}
 - 2\frac{a_2}{|c'|^5} \ip{\dot{c}'}{\dot{c}''}\,, \quad
t_3 = 2\frac{a_1-b_1}{|c'|^3} \ip{\dot{c}'}{c'}\,, \quad
t_4 = 2a_0 |c'| \,,\\
t_5 &= 2\frac{a_2}{|c'|^7} \ip{c'}{c''} \,,\quad
t_6 =2\frac{a_1}{|c'|}\,,\quad t_7=2\frac{b_1}{|c'|}\,\quad
t_8 = - 2\frac{a_2}{|c'|^5} \ip{c'}{c''} \,, \quad
t_9 = 2\frac{a_2}{|c'|^3}\,. 
\end{align*}
The varifold distance as a function of only its left argument is given by
\begin{equation*}
 F(c_1) = \langle \mu_{c_1} , \mu_{c_2} \rangle_{\on{Var}} = \iint_{S^1\x S^1} \rho(|c_1(\theta_1) - c_2(\theta_2)|^2) \gamma\left( \ip{v_1}{v_2} \right) \ud s_1 \ud s_2\,.
\end{equation*}
with the tangent and normal vectors defined by
\begin{equation*}
    v_1(\theta_1) = \frac{c_1'(\theta_1)}{|c_1'(\theta_1)|}, \quad v_2(\theta_2) = \frac{c_2'(\theta_2)}{|c_2'(\theta_2)|}
\end{equation*}
The variation of these quantities is simply
\[
D_{c_1,h}(v_1) = \langle D_s h, n \rangle n, \quad D_{c_1,h}(v_2) = 0.
\]
The derivative is given by the formula
\begin{align*}
     d F_{c_1} (h) &= \iint_{S^1\x S^1} 2\rho'(|c_1(\theta_1) - c_2(\theta_2)|^2)  \gamma\left( \ip{v_1}{v_2} \right) \ip{c_1(\theta_1) - c_2(\theta_2)}{h(\theta_1)} \ud s_1 \ud s_2 \\
     & \quad + \iint_{S^1\x S^1} \rho(|c_1(\theta_1) - c_2(\theta_2)|^2) \gamma'\left( \ip{v_1}{v_2}  \right) \ip{D_{s_1} h }{n_1} \ip{n_1}{v_2}  \ud s_1 \ud s_2 \\
     & \quad + \iint_{S^1\x S^1} \rho(|c_1(\theta_1) - c_2(\theta_2)|^2) \gamma\left( \ip{v_1}{v_2} \right) \ip{D_{s_1} h }{v_1} \ud s_1 \ud s_2\,.
\end{align*}

\end{document}